\documentclass[times,final,12pt]{article}

\usepackage{framed,multirow}
\usepackage{amssymb,latexsym,amsmath,amsthm}
\usepackage{array,booktabs,makecell,rotating}
\usepackage{caption,subcaption}
\usepackage{graphics}
\usepackage{url}
\usepackage{xcolor}
\usepackage{lineno}
\usepackage{hyperref}
\usepackage{verbatim}
\usepackage{soul}
\usepackage{algorithmicx,algorithm}
\usepackage[%
paperwidth=21.59cm,%
paperheight=27.94cm,%
margin=2in,%
left=2cm,%
right=2cm,%
top=2.5cm,%
bottom=2.5cm
]{geometry}
\usepackage{tikz}
\usepackage{subfiles}
\usepackage{booktabs}
\usepackage{xcolor}

\definecolor{newcolor}{rgb}{.8,.349,.1}
\newcommand{\red}{\color{red}}

\newcommand{\be}{\begin{eqnarray}}
\newcommand{\ee}{\end{eqnarray}}
\newcommand{\bes}{\begin{eqnarray*}}
\newcommand{\ees}{\end{eqnarray*}}
\newcommand{\beqn}{\begin{equation}}
\newcommand{\eeqn}{\end{equation}}
\newcommand{\beqns}{\begin{equation*}}
\newcommand{\eeqns}{\end{equation*}}
\newcommand{\bitem}{\begin{itemize}}
\newcommand{\eitem}{\end{itemize}}

\newcommand{\dx}{{\Delta x}}

\newcommand{\sst}{\scriptscriptstyle}

\newcommand{\beqa}{\begin{eqnarray}}
\newcommand{\eeqa}{\end{eqnarray}}

\newcommand{\bF}{{\mathbf{F}}}

\newcommand{\bQ}{{\mathbf{Q}}}

\newcommand{\half}{\frac 1 2}

\newcommand{\MajorHead}[1]{\bigskip\begin{center}{\underline{\bf #1}}
                                   \end{center}}

\renewcommand{\thebibliography}[1]{\MajorHead{References}
                  \list {[\arabic{enumi}]}{\settowidth\labelwidth{[#1]}
                         \leftmargin\labelwidth
                         \advance\leftmargin\labelsep
                         \usecounter{enumi}}
                         \def\newblock{\hskip .11em plus .33em minus -.07em}
                         \sloppy \sfcode`\.=1000\relax}

\definecolor{background}{cmyk}{1,1,0,0}

\newtheorem{theorem}{Theorem}

\newtheorem{tmm}{Theorem}

\newtheorem{example}[theorem]{Example}
\newtheorem{remark}[tmm]{Remark}

\newcommand{\ncp}{n_{\sst cp}}

\newcommand{\bff}{\mathbf{f}}

\newcommand{\bfc}{\mathbf{c}}

\newcommand{\bfd}{\mathbf{d}}

\newcommand{\intd}{\mathrm{d}}

\newcommand{\CD}{\mathrm{CD}}
\newcommand{\JS}{\mathrm{JS}}
\newcommand{\M}{\mathrm{M}}
\newcommand{\Z}{\mathrm{Z}}

\newcommand{\eps}{\varepsilon}

\newcommand{\order}{{\mathcal{O}}}

\newcommand{\e}{\mathrm e}

\newcommand{\minus}{\!-\!}
\newcommand{\cfl}{\text{CFL}}

\newcommand{\Real}{\text{Re}}
\newcommand{\Imag}{\text{Im}}

\title{
A Family of Even-Order Central-Upwind WENO Schemes with Averaged Downwind and Novel Global Smoothness Indicators
}
\author{
Jiaxi Gu\footnote{Department of Mathematics and POSTECH MINDS (Mathematical Institute for Data Science), Pohang University of Science and Technology, Pohang, Korea. E-Mail: jiaxigu@postech.ac.kr}
, \
Bao-Shan Wang\footnote{Corresponding author: School of Mathematical Sciences, Ocean University of China, Qingdao, China. E-Mail: wbs@stu.ouc.edu.cn}
, \
Wai Sun Don\footnote{Department of Mathematics, Hong Kong Baptist University, Hong Kong. E-Mail: donwaisun@outlook.com}
, \
Jae-Hun Jung\footnote{Department of Mathematics and POSTECH MINDS (Mathematical Institute for Data Science), Pohang University of Science and Technology, Pohang, Korea. E-Mail: jung153@postech.ac.kr}
}
\date{}

\begin{document}
\baselineskip=1.5pc
\maketitle

\begin{abstract}
We propose a simple yet effective local smoothness indicator for the downwind stencil in central-upwind weighted essentially non-oscillatory (WENO) schemes of even order for hyperbolic conservation laws. Starting from an odd-order upwind WENO scheme, we construct an even-number-of-points stencil by incorporating a downwind substencil whose smoothness indicator is the arithmetic mean of all local smoothness indicators. This straightforward averaging approach incorporates regularity information from the entire stencil without requiring additional tuning parameters or complex formulations. Combined with affine-invariant Z-type nonlinear weights and a carefully designed global smoothness indicator, the resulting scheme, termed WENO-ZA6 for the sixth-order case, achieves optimal convergence rates at critical points up to second order, exhibits favorable dispersion and dissipation properties as confirmed by approximate dispersion relation analysis, and provides sharp, essentially non-oscillatory resolution of discontinuities. Numerical experiments on scalar problems and the one- and two-dimensional Euler equations demonstrate that WENO-ZA6 achieves accuracy comparable to or better than existing sixth-order central-upwind schemes (WENO-CU6, WENO-S6) and the seventh-order WENO-Z7, while requiring approximately 15\%--21\% less computational time. The framework extends naturally to fourth-, eighth-, and tenth-order schemes.
\end{abstract}



\section{Introduction}
\label{sec:introduction}

This work addresses the solution of hyperbolic conservation laws in the form of
\begin{equation} \label{eq:hyperbolic}
\frac{\partial{\bQ}}{\partial{t}} + \nabla \cdot \bF(\bQ) = 0,
\end{equation} 
where $\bF(\bQ)$ is the flux function of the conservative variables $\bQ$. 
Such equations are fundamental to diverse applications, including gas dynamics, magnetohydrodynamics, shallow-water flows, and traffic modeling. 
Even for smooth initial data, solutions may develop finite-time singularities such as shocks, rarefaction waves, and contact discontinuities. 
A primary challenge for numerical methods is to achieve high-order accuracy in smooth regions while maintaining sharp, essentially non-oscillatory (ENO) resolution across discontinuities. These objectives are inherently conflicting.

The weighted essentially non-oscillatory (WENO) methodology, introduced by Liu, Osher, and Chan \cite{Liu} in the finite volume framework and extended to finite differences by Jiang and Shu \cite{Jiang}, has become one of the most successful approaches for addressing this challenge. 
The central concept involves constructing the numerical flux as a nonlinear convex combination of candidate fluxes, selected from multiple substencils with weights that adapt to local smoothness.
In smooth regions, the weights approximate their ideal linear values, thereby recovering the optimal order of accuracy associated with central approximations.
Near discontinuities, substencils containing the discontinuity receive negligible weights. The scheme then adaptively reduces to a lower-order, yet essentially non-oscillatory, reconstruction on the remaining smooth substencils.

The original WENO-JS5 scheme, developed by Jiang and Shu \cite{Jiang}, achieves fifth-order accuracy in smooth regions by forming a convex combination of three third-order candidate stencils. 
However, the classical (nonlinear) weights $\omega_k^\JS$ fail to recover the ideal linear weights with sufficient accuracy at first-order critical points, where $f'(x) = 0$ and $f''(x) \ne 0$, resulting in a reduction of accuracy to third order.
Henrick et al. \cite{Henrick} addressed this limitation by introducing mapping functions that map the classical weights $\omega_k^\JS$ to new mapped weights $\omega_k^\M$ satisfying the sufficient conditions for optimal convergence. 
Alternatively, Borges et al. \cite{Borges} introduced a global smoothness indicator $\tau$ to measure the regularity of the entire stencil and defined Z-type weights $\omega_k^\Z$, which achieve optimal order at critical points with lower computational cost than mapped weights.
Subsequent developments have further refined these ideas, including affine-invariant formulations \cite{DonLiWangWang, WangDon} that maintain adaptivity for problems with small-amplitude discontinuities.

Some studies aim to reduce the numerical dissipation inherent in purely upwind WENO schemes. By incorporating information from the downwind direction, central-upwind WENO schemes achieve higher-order accuracy on a given stencil while preserving the ENO property near discontinuities.
Hu et al. \cite{HuWangAdams} developed an adaptive central-upwind scheme (WENO-CU6) that attains sixth-order accuracy by including a downwind substencil in the reconstruction. 
Similar to the Z-type weights, the central-upwind adaptation is controlled by assigning a smoothness indicator to the downwind substencil that incorporates global regularity information. In smooth regions, the downwind substencil contributes to the reconstruction, yielding sixth-order central accuracy. Near discontinuities, the downwind contribution is suppressed, and the scheme reverts to fifth-order upwind bias. This framework was extended by Hu et al. \cite{HuAdams} for implicit large-eddy simulation and by Huang and Chen \cite{Huang}, who proposed WENO-S6 with a different construction of the global smoothness indicator. Additional sixth-order central-upwind schemes have been developed in \cite{ZhaoDuYuan,ZhaoSunMei,Wang}, each with different choices for the downwind smoothness indicator or weight formulation. Readers are referred to the review literature for an extensive recent state-of-the-art development of the WENO methodology.  For example, the article “Recent Advancements in Fluid Flow Simulation Using the WENO Scheme: A Comprehensive Review, 2024” \cite{Bozorgpour2025} gives a good introduction to the most recent developments in the field. 

Despite their success, existing central-upwind WENO schemes depend on complex downwind smoothness indicators and often require carefully tuned combinations of global and local measures. This complexity increases implementation effort and computational cost, and optimal parameter choices are often problem-dependent. This raises the question: can a simpler construction achieve comparable or superior performance?

Recent work \cite{GuChenParkJung} provided an affirmative answer for the fourth-order case by defining the downwind smoothness indicator as the arithmetic mean of all local smoothness indicators.
This simple averaging incorporates regularity information from the entire stencil. A large smoothness indicator from any substencil increases the average, thereby reducing the influence of the downwind candidate. The fourth-order scheme demonstrated excellent performance with minimal implementation complexity.

The present study extends this averaging approach to central-upwind WENO schemes of arbitrary even order, with schemes up to tenth order presented. The construction begins with a $(2r-1)$th-order upwind WENO scheme on a $(2r-1)$-point stencil, followed by the addition of one downwind substencil to form a $2r$-point stencil.
The smoothness indicator for the downwind substencil is redefined as the average of all $r+1$ local smoothness indicators, including those from the $r$ upwind substencils and the downwind substencil itself. When combined with Z-type affine-invariant weights and a global smoothness indicator based on discrete high-order derivative approximations \cite{WangDonLiWang}, this approach yields a family of even-order central-upwind WENO schemes with the following properties:
\begin{itemize}
\item \textbf{Simplicity:} The downwind smoothness indicator requires no additional tuning parameters beyond those already present in the weight formulation.

\item \textbf{Optimal accuracy:} For the sixth-order representative scheme (termed WENO-ZA6), we prove that the nonlinear weights satisfy the sufficient conditions for sixth-order convergence at critical points up to second order with $p = 1$, and up to third order with $p = 3/2$.  

\item \textbf{Favorable spectral properties:} Approximate dispersion relation analysis confirms that WENO-ZA6 has dispersion characteristics close to the underlying sixth-order linear central scheme. It also maintains strictly positive numerical dissipation across all wavenumbers with the proposed global smoothness indicator.

\item \textbf{Computational efficiency:} The simplified construction reduces computational cost by approximately 15–21\% compared to WENO-S6 and WENO-Z7, while achieving comparable or superior accuracy.

\item \textbf{Robustness:} Extensive numerical experiments on scalar equations and the Euler system in one and two dimensions demonstrate sharp shock capturing without spurious oscillations and resolving fine-scale structures efficiently.
\end{itemize}

The remainder of this paper is organized as follows:
In Section~\ref{sec:CUWENO}, we describe the general framework for central-upwind WENO schemes of arbitrary even order, including the stencil decomposition, smoothness indicators, and generalized nonlinear weights. 
Section~\ref{sec:CUWENO6} focuses on the sixth-order WENO-ZA6 scheme as a representative example, providing detailed convergence analysis at critical points and approximate dispersion relation analysis. 
Numerical results for scalar problems and the Euler equations in one and two dimensions are presented in Section~\ref{sec:numerics}, comparing WENO-ZA6 with WENO-CU6 \cite{HuWangAdams}, WENO-S6 \cite{Huang}, and WENO-Z7 \cite{DonBorges}. 
Concluding remarks are given in Section~\ref{sec:conclusion}. 
For completeness, Appendix~\ref{sec:appendix} provides explicit formulas for the fourth-, eighth-, and tenth-order schemes.

\section{Even-order central-upwind WENO scheme}
\label{sec:CUWENO}

We consider the one-dimensional scalar hyperbolic conservation law,
\begin{equation}\label{eq:hyperbolic_1D_scalar}
\frac{\partial{Q}}{\partial{t}}+\frac{\partial{f(Q)}}{\partial{x}} = 0, \qquad x \in \Omega=[a,b], \quad t > 0,
\end{equation}
subject to the initial condition $Q(x,0) = Q_0(x)$ and appropriate boundary conditions.
In the following discussion, we omit the dependence on time $t$ and focus on the spatial discretization.
We apply a uniform mesh to the computational domain $\Omega=[a,b]$ with cells 
$$ I_i = [x_{i-\half}, x_{i+\half}], \quad x_{i+\half} = a + i\dx, \quad i = 1, \dots, N, $$ 
where $\dx = x_{i+\half}-x_{i-\half} = (b-a)/N$ is the cell size, $x_i=\half(x_{i-\half}+x_{i+\half})$ is the cell center, and $N$ is the number of cells.
Define the auxiliary function $h(x)$ implicitly by 
\begin{equation}\label{eq:auxiliary_func}
f \left( Q(x) \right) = \frac{1}{\dx} \int_{x-\dx/2}^{x+\dx/2} h(\xi) \, \intd \xi.
\end{equation}
Differentiating both sides with respect to $x$, we obtain
\begin{equation} \label{eq:partial_f}
\frac{\partial f}{\partial x} = \frac{h(x+\Delta x/2) - h(x-\Delta x/2)}{\dx}.
\end{equation}
With the evaluation of \eqref{eq:partial_f} at $x=x_i$, \eqref{eq:hyperbolic_1D_scalar} becomes 
\begin{equation} \label{eq:hyperbolic_1D_scalar_h}
\frac{\intd Q_i(t)}{\intd t} = - \frac{h_{i+\half} - h_{i-\half}}{\dx},
\end{equation} 
with $Q_i(t)=Q(x_i,t)$ and $h_{i \pm \half} = h(x_{i \pm \half})$.
Our goal is to reconstruct the fluxes $h_{i \pm \half}$ at the cell interfaces so that the resulting approximation achieves high-order accuracy in smooth regions while avoiding spurious oscillations near discontinuities.

To approximate the flux $h_{i+\half}$, we first use the Lax-Friedrichs flux splitting:
\begin{equation}\label{eq:LF}
f^{\pm}(Q) = \half \left( f(Q) \pm \alpha Q \right), \quad \alpha = \max_{Q} |f'(Q)|,
\end{equation}
and then construct the positive and negative numerical fluxes $\hat{f}^{\pm}_{i+\half}$ at the cell interface $x_{i+\half}$ with 
$$ \hat{f}_{i+\half} = \hat{f}^+_{i+\half} + \hat{f}^-_{i+\half}. $$
It suffices to consider the reconstruction of the positive numerical flux $\hat{f}^+_{i+\half}$, since the reconstruction of the negative numerical flux $\hat{f}^-_{i+\half}$ can be obtained by mirror symmetry with respect to the cell interface $x_{i+\half}$. 
Thus, the superscript $+$ will be omitted for simplicity whenever no confusion arises.

\begin{figure}[h!]
\vspace{-0.15in}
\centering
\begin{tikzpicture}
 \draw (-7,2.3) -- (7,2.3);
 \filldraw[red] (-7,2.3) circle (2pt);
 \filldraw[red] (-5,2.3) circle (2pt);
 \filldraw[red] (-3,2.3) circle (2pt);
 \filldraw[red] (-1,2.3) circle (2pt);
 \filldraw[red] (1,2.3)  circle (2pt);
 \filldraw[red] (3,2.3)  circle (2pt);
 \filldraw[red] (5,2.3)  circle (2pt);
 \filldraw[red] (7,2.3)  circle (2pt);

 \node[rectangle, inner sep=2pt, minimum size=2.5pt, draw=blue] (12) at (0,2.3) {};
 \draw [->,thick] (0,3.03)--(12);
 \draw [blue,->,thick] (-1,2.73)--(0,2.73);
 \node at (-0.55,3) {\scriptsize $\displaystyle \textcolor{blue}{\hat{f}^+_{i+\half}}$};

 \node at (-7+0.12,2) {$x_{i-3}$};
 \node at (-5+0.12,2) {$x_{i-2}$};
 \node at (-3+0.12,2) {$x_{i-1}$};
 \node at (-1-0.02,2) {$x_{i}$};
 \node at (1+0.13,2)  {$x_{i+1}$};
 \node at (3+0.12,2)  {$x_{i+2}$};
 \node at (5+0.12,2)  {$x_{i+3}$};
 \node at (7+0.12,2)  {$x_{i+4}$};
 
 \draw (-7,1.56) -- (7,1.56);
 \node at (-7.28,1.62) {\footnotesize $S^8$}; 
 \filldraw[green] (-7,1.56) circle (2pt);
 \filldraw[green] (-5,1.56) circle (2pt);
 \filldraw[green] (-3,1.56) circle (2pt);
 \filldraw[green] (-1,1.56) circle (2pt);
 \filldraw[green] (1,1.56)  circle (2pt);
 \filldraw[green] (3,1.56)  circle (2pt);
 \filldraw[green] (5,1.56)  circle (2pt);
 \filldraw[green] (7,1.56)  circle (2pt);
	
 \draw (-7,1.2) -- (5,1.2);
 \node at (-7.28,1.25) {\footnotesize $S^7$};
 \filldraw[blue] (-7,1.2) circle (2pt);
 \filldraw[blue] (-5,1.2) circle (2pt);
 \filldraw[blue] (-3,1.2) circle (2pt);
 \filldraw[blue] (-1,1.2) circle (2pt);
 \filldraw[blue] (1,1.2)  circle (2pt);
 \filldraw[blue] (3,1.2)  circle (2pt);
 \filldraw[blue] (5,1.2)  circle (2pt);
	
 \draw (-7,0.8) -- (-1,0.8);
 \node at (-7.28,0.8)   {\footnotesize $S_0$};
 \node at (-1+0.35,0.8) {\footnotesize $\omega_0$};
 \filldraw[black] (-7,0.8) circle (2pt);
 \filldraw[black] (-5,0.8) circle (2pt);
 \filldraw[black] (-3,0.8) circle (2pt);
 \filldraw[black] (-1,0.8) circle (2pt);
	
 \draw (-5,0.4) -- (1,0.4);
 \node at (-5.28,0.4) {\footnotesize $S_1$};
 \node at (1.35,0.4)  {\footnotesize $\omega_1$};
 \filldraw[black] (-5,0.4) circle (2pt);
 \filldraw[black] (-3,0.4) circle (2pt);
 \filldraw[black] (-1,0.4) circle (2pt);
 \filldraw[black] (1,0.4)  circle (2pt);
	
 \draw (-3,0) -- (3,0);
 \node at (-3.28,0) {\footnotesize $S_2$};
 \node at (3.35,0) {\footnotesize $\omega_2$};
 \filldraw[black] (-3,0) circle (2pt);
 \filldraw[black] (-1,0) circle (2pt);
 \filldraw[black] (1,0) circle (2pt);
 \filldraw[black] (3,0) circle (2pt);

 \draw (-1,-0.4) -- (5,-0.4);
 \node at (-1.28,-0.4) {\footnotesize $S_3$};
 \node at (5.35,-0.4) {\footnotesize $\omega_3$};
 \filldraw[black] (-1,-0.4) circle (2pt);
 \filldraw[black] (1,-0.4) circle (2pt);
 \filldraw[black] (3,-0.4) circle (2pt);
 \filldraw[black] (5,-0.4) circle (2pt);
 
 \draw[dashed] (1.1,-0.8) -- (2.9,-0.8);
 \draw[dashed] (3.1,-0.8) -- (4.9,-0.8);
 \draw[dashed] (5.1,-0.8) -- (6.9,-0.8);
 \node at (1-0.28,-0.8) {\footnotesize $S_4$};
 \node at (7.35,-0.8) {\footnotesize $\omega_4$};
 \draw (1,-0.8) circle (2pt);
 \draw (3,-0.8) circle (2pt);
 \draw (5,-0.8) circle (2pt);
 \draw (7,-0.8) circle (2pt);
\end{tikzpicture}
\caption{The stencil $S^8$ together with its five $4$-point substencils $\{S_0,S_1,S_2,S_3\}$ (upwind) and $S_4$ (downwind) used in the reconstruction of $\hat f^+_{i+\half}$ for the eighth-order central-upwind WENO scheme. The corresponding configuration for $\hat f^-_{i+\half}$ is obtained by mirror symmetry.}
\label{fig:stencil8}
\end{figure}
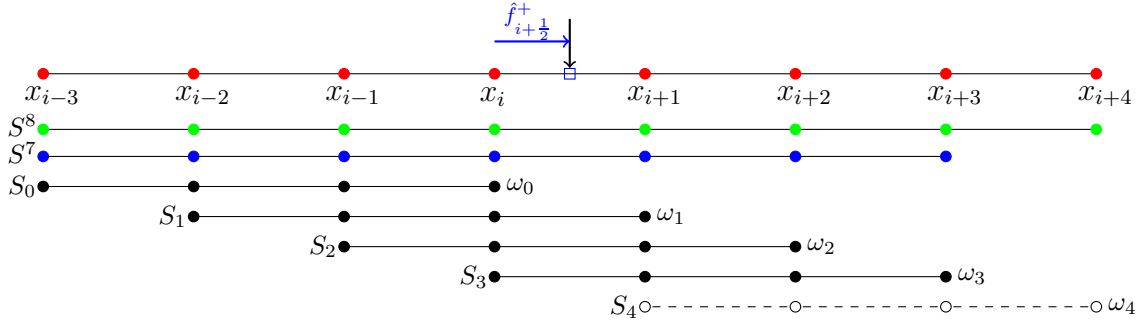

For the $2r$th-order central reconstruction of the numerical flux $\hat{f}_{i+\half}$, we consider the $2r$-point stencil,
\begin{equation}\label{eq:stencil}
S^{2r} = \{ x_{i-r+1}, x_{i-r+2}, \dots, x_{i+r-1}, x_{i+r} \}.
\end{equation}
For comparison, the standard upwind reconstruction utilizes the $(2r \minus 1)$-point stencil,
$$ S^{2r-1} = S^{2r} \setminus \{ x_{i+r} \}. $$
On the stencil $S^{2r}$, we define the polynomial $p(x)$ of degree at most $2r-1$ to approximate the auxiliary function $h(x)$:
\begin{equation}\label{eq:poly_h}
p(x) = \sum_{n=0}^{2r-1} a_n 
\left( \frac{x-x_i}{\dx} \right)^n.
\end{equation}
The coefficients $a_n$ are determined by 
\begin{equation}\label{eq:poly_h_coeffs}
f_{i+\ell} 
= \frac{1}{\dx} 
\int_{I_{i+\ell}} 
p(\xi) \, \intd \xi
= \sum_{n=0}^{2r-1} a_n \frac{1}{n+1} 
\left[ \left( \ell+\half \right)^{n+1} - \left( \ell-\half \right)^{n+1} \right], \quad   \ell = -r+1, \dots, r,
\end{equation}
where $f_{i+\ell} = f(Q_{i+\ell})$.
Evaluating $p(x)$ at $x = x_{i+\half}$ gives the central difference (CD) numerical flux
\begin{equation} \label{eq:numerical_flux_FD_plus}
\hat{f}^{\CD}_{i+\half} = p(x_{i+\half}) = \langle \bfc_{2r}, \bff^{2r}_+ \rangle, \quad \bff^{2r}_+ = (f_{i-r+1}, \dots, f_{i+r}),
\end{equation}
where $\bfc_{2r}$ denotes the coefficient vector associated with the stencil $S^{2r}$.
Next, we decompose the stencil $S^{2r}$ into $r+1$ substencils, consisting of the upwind substencils $\{S_0, \dots, S_{r-1}\}$ and the downwind substencil $S_r$, where
$$ S_k = \{ x_{i-r+1+k}, \dots, x_{i+k} \}, \quad k=0, \dots, r. $$
For $r=4$, the decomposition of $S^8$ into five substencils is illustrated in Fig. \ref{fig:stencil8}.
On each substencil $S_k$, we define a polynomial $p_k(x)$ of degree at most $r-1$.
Evaluating $p_k(x)$ at $x=x_{i+\half}$ gives the $r$th-order candidate numerical fluxes
\begin{equation}\label{eq:numerical_subflux}
\hat{f}^k_{i+\half} = p_k(x_{i+\half}) = \langle \bfc_k, \bff_k \rangle, \quad \bff_k = (f_{i-r+1+k}, \dots, f_{i+k}),
\end{equation}
where $\bfc_k$ denotes the coefficient vector on $S_k$.
A linear combination of these candidate numerical fluxes $\hat{f}^k_{i+\half}$ recovers the central difference numerical flux $\hat{f}^{\CD}_{i+\half}$ in \eqref{eq:numerical_flux_FD_plus}. 
More precisely,
\begin{equation}\label{eq:numerical_flux_FD_linear_comb}
\langle \bfc_{2r}, \bff^{2r}_+ \rangle = \langle \bfd, \hat{\bff} \rangle, \quad \hat{\bff} = (\hat{f}^0_{i+\half}, \dots, \hat{f}^r_{i+\half}),
\end{equation}
where $\hat{\bff}$ is the vector of candidate numerical fluxes and $\bfd = (d_0, \dots, d_r)$ is the vector of ideal weights obtained by solving the system of equations \eqref{eq:numerical_flux_FD_linear_comb}, with the explicit form
\begin{equation}
    d_k = \binom{r}{k}^2\bigg/\binom{2r}{r},\quad k = 0,1,\dots,r,\quad \binom{m}{n} = \frac{m!}{n! \, (m-n)!}.
\end{equation}

The WENO numerical flux $\hat{f}_{i+\half}$ takes a convex combination of $r+1$ candidate numerical fluxes $\hat{f}^k_{i+\half}$,
\begin{equation}\label{eq:numerical_flux_WENO}
\hat{f}_{i+\half} = \langle \boldsymbol{\omega}, \hat{\bff} \rangle = \sum_{k=0}^r \omega_k \, \hat{f}^k_{i+\half},
\end{equation}
where the nonlinear weights $\boldsymbol{\omega} = (\omega_0, \dots, \omega_r)$ satisfies $\omega_k \ge 0$ and $\sum \omega_k = 1$.
In the central-upwind framework, the downwind candidate $\hat{f}^r_{i+\half}$, constructed on $S_r$, is not used merely as part of a high-order polynomial reconstruction, but is instead treated as an independent directional sensor. 
This allows the scheme to switch between a highly accurate central approximation, which fully uses $S^{2r}$ in smooth regions, and a robust upwind approximation, which reverts to $S^{2r-1}$ near discontinuities.

\subsection{Generalized WENO weights}
\label{sec:weights}

We now define the smoothness indicators and the generalized nonlinear weights $\omega_k$ in \eqref{eq:numerical_flux_WENO} for the central-upwind reconstruction. 
Specifically, we introduce local smoothness indicators on the substencils $S_k$ and a global smoothness indicator on the full stencil $S^{2r}$.

The classical (local) smoothness indicators $\beta_k$ \cite{Jiang} are constructed by measuring the scaled $L^2$ norms of the derivatives of the local reconstructed polynomials $p_k(x)$ over the substencils $S_k$.
Using \eqref{eq:poly_h}, and taking $\dx=1$ and $ x_i=0$ for convenience, the smoothness indicators $\beta_k, \; k \in \{ 0, \dots, r-1, d\}$ can be expressed in terms of the coefficients $a_n$ as
$$
\beta_k = \sum_{\ell=1}^{r-1} \dx^{2\ell-1} \int_{I_i} \left( \frac{\intd^\ell p_k(x)}{\intd x^\ell} \right )^2 \intd x 
        = \sum_{\ell=1}^{r-1} \int_{-\half}^{\half} \left( \sum_{n=\ell}^{r-1} B_n^{\ell} \xi^{n-\ell} \right)^2 \intd \xi,   
$$
and hence, by setting $\kappa= i+j-2\ell$,
\begin{equation}\label{eq:SI}
\beta_k = \sum_{\ell=1}^{r-1} \sum_{\substack{i+j=2\ell,~\textbf{even}, \\ ~\ell \leq i,j \leq r-1}}^{2(r-1)} C_{i,j}^{\ell} B_{i}^{\ell} B_{j}^{\ell}, \quad \text{with} \quad B_n^{\ell} = \frac{n!}{(n-\ell)!}\, a_n, \; C_{i,j}^{\ell} = \frac {1} {\kappa+1}2^{-\kappa}.
\end{equation}
Here, $\beta_d$ denotes the downwind smoothness indicator obtained from the polynomial $p_r(x)$ on the downwind substencil $S_r$.
The smoothness indicator $\beta_r$ is then defined as the arithmetic mean of all smoothness indicators \cite{GuChenParkJung},
\begin{equation}\label{eq:SI_downwind} 
\beta_r = \frac{1}{r+1} \left( \sum_{k=0}^{r-1} \beta_k + \beta_d \right). 
\end{equation}
The advantage of the new definition is that it enables the downwind smoothness indicator to incorporate regularity information from the entire stencil \(S^{2r}\), rather than relying solely on a local substencil.

The global smoothness indicator on \(S^{2r}\) is defined as the square of a discrete \((2r-1)\)-th derivative measure~\cite{WangDonLiWang}:
\begin{equation}\label{eq:GSI}
\tau = \left\langle \mathbf{c}_\tau, \mathbf{f}^{2r} \right\rangle^2,
\end{equation}
where the coefficient vector \(\mathbf{c}_\tau\) is chosen so that 
\(\left\langle \mathbf{c}_\tau, \mathbf{f}^{2r} \right\rangle\) 
provides a consistent approximation of the \((2r-1)\)-th derivative at \(x_i\).
With this construction, the leading-order behavior in smooth regions satisfies
$
\tau = \mathcal{O}\!\left(\Delta x^{4r-2}\right),
$
provided that \(f_i^{(2r-1)} \neq 0\); that is, when the first nonvanishing derivative at \(x_i\) is of odd order.

\begin{remark}\rm
An alternative global smoothness indicator may be defined by combining discrete $(2r \minus 1)$th- and $(2r \minus 2)$th-derivative measures,
\begin{equation}\label{eq:GSI_alt}
\tau' = \frac{13}{12} \langle \bfc_\tau, \bff^{2r} \rangle^2 + \langle \bfc_\tau', \bff^{2r} \rangle^2,
\end{equation}
whose leading-order behavior in smooth regions is $\tau=\order(\dx^{4r-4})$ if $f^{(2r-2)}_i \ne 0$, that is, when the first nonvanishing derivative is of even order. 
The first term measures the smoothness of the full stencil $S^{2r}$, while the second term corresponds to the smoothness of the stencil $S^{2r-1}$.
As shown in Section \ref{sec:ADR}, this choice tends to introduce antidiffusion, so we do not use it in the present work.
\end{remark}

With the local and global smoothness indicators in place, we are ready to define the nonlinear weights $\omega_k$ for the numerical flux \eqref{eq:numerical_flux_WENO}.
The nonlinear weights $\omega_k,~k=0,\dots,r$ take the form
\begin{equation}\label{eq:weights_general}
\omega_k = \frac{\alpha_k}{\sum_{s=0}^r \alpha_s}, 
\end{equation}
where the unnormalized weights $\alpha_k$ are defined compactly as
\begin{equation}\label{eq:weights_unnormalized}
\alpha_k = d_k \left( L + \Gamma_k \right)^q.
\end{equation}
Here, $d_k$ are the ideal linear weights in \eqref{eq:numerical_flux_FD_linear_comb}.
The term inside the parentheses consists of a constant linear component $L$ and a nonlinear component $\Gamma_k$. 
The parameter $q \geq 1$ can be used to enhance the separation between the linear and nonlinear components.
The components are defined as follows:
\begin{enumerate}
\item \textbf{Linear Component ($L$):} This is a constant, typically set to unity ($L=1$), representing the baseline contribution of the ideal weights in smooth regions and independent of the substencil $S_k$.

\item \textbf{Nonlinear Component ($\Gamma_k$):} This term measures the ratio of global to local smoothness. 
It reduces the contribution of the local polynomial on a nonsmooth substencil $S_k$ in the reconstruction procedure,
\begin{equation}\label{eq:weight_nonlinear}
\Gamma_k = \left( \frac{\tau}{\beta_k + \eps \cdot \mu[f]^2} \right)^p,
\end{equation}
where $\eps > 0$ is the sensitivity parameter and $p \geq 1$ is the power parameter. 
For the affine-invariant (Ai) property, the descaler is defined by
\begin{equation}\label{eq:aimu}
\mu[f] = \frac{1}{2r} \sum_{x_j\in S^{2r}} \left| f_j-\frac{1}{2r} \sum_{x_j\in S^{2r}} f_j \right| + \mu_0,
\end{equation}
with $\mu_0>0$ a tiny constant (e.g., $10^{-40}$ in double precision and $10^{-100}$ in quadruple precision) to prevent division by zero. The reader is referred to \cite{DonLiWangWang,WangDon} for a detailed discussion about the Ai property, which is critical for preventing the potential loss of adaptivity; thus, losing the important ENO property in case of small amplitude discontinuities/shocks and generating small Gibbs oscillations.
\end{enumerate}
This general form, $\alpha_k = d_k (L + \Gamma_k)^q$, is powerful because it encapsulates distinct state-of-the-art families of both classical and affine-invariant WENO weights through appropriate parameter choices. 
Table~\ref{tab:SpecialCases} details representative parameter choices corresponding to several WENO nonlinear weights. 
\begin{table}[h!]
\centering
\caption{Representative parameter choices for WENO unnormalized weights $\alpha_k = d_k (L + \Gamma_k)^q$.}
\label{tab:SpecialCases}
\renewcommand{\arraystretch}{1.3}
\begin{tabular}{llllll}
\toprule
WENO weight & Type & $\mu[f]$ & $L$ & $p$ & $q$ \\
\midrule
Classical JS \cite{Jiang} & JS & $1$ & $1$ & $2$      & $1$ \\
Classical Z \cite{Borges} & Z  & $1$ & $1$ & $(1, 2)$ & $1$ \\
\midrule
Ai-JS \cite{WangDon}      & Ai-JS & \eqref{eq:aimu} & $0$ & $2$      & $1$ \\
Ai-Z  \cite{WangDon}       & Ai-Z  & \eqref{eq:aimu} & $1$ & $(1, 2)$ & $1$ \\
\bottomrule
\end{tabular}
\end{table}
Throughout the following discussion, we fix $L=1$ and $q=1$ for simplicity.  

\subsection{Weights behavior around discontinuity}\label{sec:weights_discontinuity}
According to the WENO strategy, when a discontinuity is contained within one or more substencils, the corresponding nonlinear weights should vanish asymptotically.

Let $D$ and $S$ denote the indices of a discontinuous substencil and a smooth substencil, respectively. 
From the smoothness indicators, we have the asymptotic estimates 
$$
\beta_D = \order(1), \qquad \beta_S = \order(\dx^2) \ \text{if } f'_i \ne 0, \qquad \tau = \order(1). 
$$
Hence, the ratio of the corresponding weights, with the setting $\eps = 0$, satisfies
\begin{align}
\frac{d_S}{d_D} \frac{\omega_D}{\omega_S} 
= \frac{1 + \left( \frac{\tau}{\beta_D} \right)^p}{1 + \left( \frac{\tau}{\beta_S} \right)^p}
= \order(\dx^{2p}) \to 0 \quad \text{as } \dx \to 0.
\end{align}
This confirms that the weight of the discontinuous substencil becomes asymptotically negligible compared to that of the smooth substencil. 

Next, we investigate the weight ratio relative to the downwind substencil $S_r$. 
Let $D$ denote the index of the discontinuous substencil with the largest smoothness indicator.
Since the smoothness indicator $\beta_r$ for the downwind substencil is defined from the average of all local smoothness indicators, it is also affected by the discontinuity, and we write $\beta_r = c \, \beta_D$ for $0 < c \leq 1$. 
Defining $\Theta = \left( \tfrac{\tau}{\beta_D} \right)^p = \order(1)$, we obtain
\begin{align}
\frac{d_r}{d_D} \frac{\omega_D}{\omega_r} 
= \frac{1 + \left( \frac{\tau}{\beta_D} \right)^p}{1 + \left( \frac{\tau}{\beta_r} \right)^p}
= \frac{1 + \Theta}{1 + c^{-p} \Theta} 
\in (c^p, 1].
\end{align}
This ratio is monotonically decreasing in $\Theta$.
Thus, the weight of the discontinuous substencil $S_D$ remains smaller than that of $S_r$, with the separation between them controlled by $c^p$. 
This allows the scheme to distinguish between substencils even when several of them are affected by the discontinuity.

In summary, discontinuous substencils are suppressed at a rate determined by $p$, whereas smooth substencils retain weights of order $\order(1)$.
The downwind candidate $k=r$ is treated in the same manner through $\beta_r$; near a discontinuity, this construction drives $\omega_r$ asymptotically zero, and the scheme reverts to an upwind bias.
Moreover, a larger value of $p$ produces a stronger separation between discontinuous and smooth substencils, leading to sharper resolution of discontinuities but potentially reduced accuracy for smooth waves at high frequencies. 
By contrast, a smaller value of $p$ yields a milder weighting, which is often more robust for smoothly varying turbulent structures, although it tends to smear discontinuities more noticeably.

\section{Sixth-order central-upwind WENO scheme}
\label{sec:CUWENO6}

In this section, we concentrate on the sixth-order case ($r=3$) as a representative example of the central-upwind WENO scheme.
This case illustrates the main ingredients of the proposed construction and will be used for numerical experiments in the following section. 
The fourth-order case $(r=2)$ can be found in \cite{GuChenParkJung}, while the cases $r=4$ and $r=5$ are provided in the Appendix.

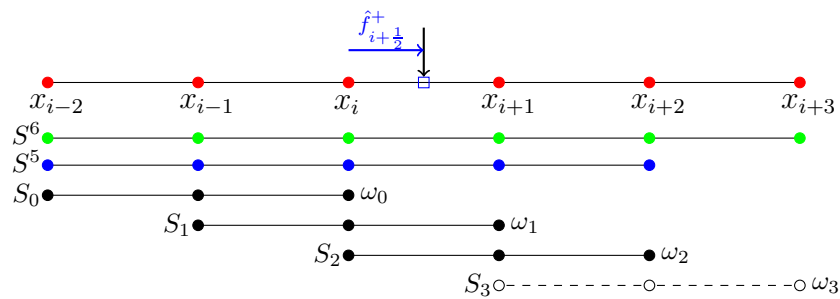
\begin{figure}[h!]
\vspace{-0.15in}
\centering
\begin{tikzpicture}
 \draw (-5,2.3) -- (5,2.3);
 \filldraw[red] (-5,2.3) circle (2pt);
 \filldraw[red] (-3,2.3) circle (2pt);
 \filldraw[red] (-1,2.3) circle (2pt);
 \filldraw[red] (1,2.3)  circle (2pt);
 \filldraw[red] (3,2.3)  circle (2pt);
 \filldraw[red] (5,2.3)  circle (2pt);

 \node[rectangle, inner sep=2pt, minimum size=2.5pt, draw=blue] (12) at (0,2.3) {};
 \draw [->,thick] (0,3.03)--(12);
 \draw [blue,->,thick] (-1,2.73)--(0,2.73);
 \node at (-0.55,3) {\scriptsize $\displaystyle \textcolor{blue}{\hat{f}^+_{i+\half}}$};

 \node at (-5+0.12,2) {$x_{i-2}$};
 \node at (-3+0.12,2) {$x_{i-1}$};
 \node at (-1-0.02,2) {$x_{i}$};
 \node at (1+0.13,2)  {$x_{i+1}$};
 \node at (3+0.12,2)  {$x_{i+2}$};
 \node at (5+0.12,2)  {$x_{i+3}$};
 
 \draw (-5,1.56) -- (5,1.56);
 \node at (-5.28,1.62) {\footnotesize $S^6$}; 
 \filldraw[green] (-5,1.56) circle (2pt);
 \filldraw[green] (-3,1.56) circle (2pt);
 \filldraw[green] (-1,1.56) circle (2pt);
 \filldraw[green] (1,1.56)  circle (2pt);
 \filldraw[green] (3,1.56)  circle (2pt);
 \filldraw[green] (5,1.56)  circle (2pt);
	
 \draw (-5,1.2) -- (3,1.2);
 \node at (-5.28,1.25) {\footnotesize $S^5$};
 \filldraw[blue] (-5,1.2) circle (2pt);
 \filldraw[blue] (-3,1.2) circle (2pt);
 \filldraw[blue] (-1,1.2) circle (2pt);
 \filldraw[blue] (1,1.2)  circle (2pt);
 \filldraw[blue] (3,1.2)  circle (2pt);
	
 \draw (-5,0.8) -- (-1,0.8);
 \node at (-5.28,0.8)   {\footnotesize $S_0$};
 \node at (-1+0.35,0.8) {\footnotesize $\omega_0$};
 \filldraw[black] (-5,0.8) circle (2pt);
 \filldraw[black] (-3,0.8) circle (2pt);
 \filldraw[black] (-1,0.8) circle (2pt);
	
 \draw (-3,0.4) -- (1,0.4);
 \node at (-3.28,0.4) {\footnotesize $S_1$};
 \node at (1.35,0.4)  {\footnotesize $\omega_1$};
 \filldraw[black] (-3,0.4) circle (2pt);
 \filldraw[black] (-1,0.4) circle (2pt);
 \filldraw[black] (1,0.4)  circle (2pt);
	
 \draw (-1,0) -- (3,0);
 \node at (-1.28,0) {\footnotesize $S_2$};
 \node at (3.35,0) {\footnotesize $\omega_2$};
 \filldraw[black] (-1,0) circle (2pt);
 \filldraw[black] (1,0) circle (2pt);
 \filldraw[black] (3,0) circle (2pt);
 
 \draw[dashed] (1.1,-0.4) -- (2.9,-0.4);
 \draw[dashed] (3.1,-0.4) -- (4.9,-0.4);
 \node at (1-0.28,-0.4) {\footnotesize $S_3$};
 \node at (5.35,-0.4) {\footnotesize $\omega_3$};
 \draw (1,-0.4) circle (2pt);
 \draw (3,-0.4) circle (2pt);
 \draw (5,-0.4) circle (2pt);
\end{tikzpicture}
\caption{The stencil $S^6$ together with its four $3$-point substencils $\{S_0,S_1,S_2\}$ (upwind) and $S_3$ (downwind) used in the reconstruction of $\hat f^+_{i+\half}$ for the sixth-order central-upwind WENO scheme. The corresponding configuration for $\hat f^-_{i+\half}$ is obtained by mirror symmetry.}
\label{fig:stencil6}
\end{figure}

We consider the $6$-point stencil $S^6$, which consists of three upwind substencils $\{ S_0, S_1, S_2 \}$ and one downwind substencil $S_3$, as shown in Fig. \ref{fig:stencil6}.
For the stencil $S^6$, the coefficient vector in \eqref{eq:numerical_flux_FD_plus} is
\begin{equation}
\bfc_6 = \tfrac{1}{60}(1,-8,37,37,-8,1).
\end{equation}
The coefficient vectors $\bfc_k$ associated with the $3$-point substencils $S_k$ are
\begin{equation}\label{eq:subflux_S6_coeffs} \textstyle
\bfc_0 = \tfrac{1}{6}(2, -7, 11), \quad
\bfc_1 = \tfrac{1}{6}(-1, 5, 2), \quad
\bfc_2 = \tfrac{1}{6}(2, 5, -1), \quad
\bfc_3 = \frac{1}{6}(11, -7, 2).
\end{equation}
The ideal weight vector in \eqref{eq:numerical_flux_FD_linear_comb} is
\begin{equation}
\bfd = \tfrac{1}{20}(1,9,9,1).
\end{equation}
Using \eqref{eq:SI}, the smoothness indicators $\beta_k$ can be written in the unified form
\begin{equation}\label{eq:SI_S6}
\beta_k
= \langle \bff_{k}, \mathbf{B}_k \bff_{k} \rangle, 
\qquad k \in \{ 0,1,2,d \},
\end{equation}
where the corresponding smoothness matrices are 
\begin{eqnarray}
\mathbf{B}_0 = \frac{1}{6}
\left[
\begin{array}{rrr}
   8 & -19 &  11 \\
 -19 &  50 & -31 \\
  11 & -31 &  20
\end{array}
\right],
\;
\mathbf{B}_1 = \frac{1}{6}
\left[
\begin{array}{rrr}
   8 & -13 &   5 \\
 -13 &  26 & -13 \\
   5 & -13 &   8
\end{array}
\right],
\;
\mathbf{B}_d = \frac{1}{6}
\left[
\begin{array}{rrr}
  44 & -73 &  29 \\
 -73 & 122 & -49 \\
  29 & -49 &  20
\end{array}
\right].
\end{eqnarray}
and $\mathbf{B}_{2} = \mathrm{rev}(\mathbf{B}_0)$, where $\mathrm{rev}(\cdot)$ denotes the reversal operator. 
The smoothness indicator for the downwind substencil $S_3$ is then defined by 
$$ \beta_3= \frac{1}{4} \left( \sum_{k=0}^2 \beta_k + \beta_d \right). $$
The global smoothness indicator $\tau$ \eqref{eq:GSI} employs the coefficient vector $\bfc_\tau = (-1,5,-10,10,-5,1)$, with $\bfc_\tau' = (1, -4, 6, -4, 1, 0)$ for the alternative version \eqref{eq:GSI_alt}.
In this paper, we employ the Z-type nonlinear weights without Ai property, and refer to the resulting scheme as WENO-ZA6.

\subsection{Convergence analysis at critical points}
For the numerical flux reconstruction at $x_{i \pm \half}$ to achieve sixth-order accuracy, a sufficient condition for the deviation of the nonlinear weights from their linear counterparts is
\begin{equation}\label{eq:WENO6_condition}
\omega_k - d_k = \order(\Delta x^4), \quad k=0,1,2,3,
\end{equation}
which ensures that the convex combination in \eqref{eq:numerical_flux_WENO} does not degrade the sixth-order accuracy of the underlying linear scheme (see \cite{DonBorges,WangDonLiWang} for details).
Defining $\bar \Gamma = \sum_j d_j \Gamma_j$ and expanding the weights in smooth regions with $\Gamma_k \ll 1$, we obtain
$$
\omega_k - d_k 
= d_k \left(\frac{1+\Gamma_k}{1+\bar \Gamma} - 1\right)
= d_k \bigl(\Gamma_k - \bar \Gamma \bigr)
+ \order(\bar \Gamma ^2), 
$$
Hence, a sufficient condition \eqref{eq:WENO6_condition} reduces to 
\begin{equation}\label{eq:Gamma_condition}
\Gamma_k = \order(\Delta x^4).
\end{equation}
The behavior of $\Gamma_k$ depends critically on the scaling of both the local and global smoothness indicators at critical points.

We assume $\eps=0$ and let $x_i$ be a critical point of order $\ncp$, (i.e., $f_i^{(1)}=\cdots=f_i^{(\ncp)}=0, \ f_i^{(\ncp+1)}\ne 0$).
It is well-known \cite{Henrick,DonBorges} that the smoothness indicators scale as
$
\beta_k = \order\!\left( \Delta x^{2(\ncp+1)} \right).
$
For the class of global smoothness indicators designed \eqref{eq:GSI}, one has
$
\tau = \order\!\left( \Delta x^{\max\{2(\ncp+1),\,10\}} \right).
$
Consequently, the sufficient condition \eqref{eq:Gamma_condition} yields
\begin{equation}
\Gamma_k
\propto
\left( \frac{\tau}{\beta_k} \right)^p
=
\left\{
\begin{array}{ll}
\order\!\left( \Delta x^{\,2p(4-\ncp)} \right),
& \ncp \le 3, \\[6pt]
\order(1),
& \ncp \ge 4.
\end{array}
\right.
\end{equation}
which leads to the explicit lower bound on the power parameter,
\begin{equation}
p \ge p^\ast := \frac{2}{4-\ncp},
\end{equation}
where $p^\ast$ denotes the \emph{critical power parameter} required to prevent order degeneration at a critical point of order $\ncp$.
For $\ncp \leq 3$, choosing $p \geq p^\ast$ guarantees that the nonlinear weights recover the optimal linear weights rapidly enough to retain sixth-order accuracy. 
Then, increasing $p$ further enhances recovery of linear weights, though at the cost of reduced sensitivity near smooth extrema.
For higher-order critical points ($\ncp \geq 4$), however, $\Gamma_k=\order(1)$ independently of $p$.
As a result, the sufficient condition \eqref{eq:Gamma_condition} cannot be satisfied, and a loss of the optimal sixth-order accuracy is inevitable. 
The resulting convergence behavior is summarized in Table~\ref{tab:WENO6_condition_ZA_p}.
\begin{table}[htb!]
\centering
\begin{tabular}{cccc}
\toprule
$\ncp$ & $p^*$ & $\lceil{p^*\rceil}$ &  $\theta \left( \frac{1}{\dx}\left(\hat f_{i+\half}-\hat f_{i-\half}\right)\right)$ \\
\midrule
1 & $\frac{2}{3}$ & 1 & 6 \\
2 & $1$           & 1 & 6 \\
3 & $\frac{3}{2}$ & 2 & 6 \\
4 & --            & --                        & 4 \\
5 & --            & --                       & 5 \\
\bottomrule
\end{tabular}
\caption{Sufficient values of $p$ for nonlinear weights to satisfy the sixth-order condition \eqref{eq:WENO6_condition} at a critical point of order $\ncp$. $\theta(g)$ denotes the leading order of accuracy for $g(\dx)$ as in $\order{(\dx^{\theta(g)})}$.}
\label{tab:WENO6_condition_ZA_p}
\end{table}
A possible remedy is to introduce a variable sensitivity parameter $\eps = \order\!\left( \Delta x^{{2(\ncp-1)}/p} \right)$ that restores the optimal order of accuracy in the presence of high-order critical points \cite{DonBorges}.
Its potential drawback is that, on coarse meshes, the loss of adaptivity (e.g., $\beta_k \ll \eps$) may lead to increased Gibbs oscillations near discontinuities. 

To verify the theoretical results above, we consider the test function
\begin{equation}\label{eq:func_convergence_critical_pt}
f(x) = x^{\ncp+1}e^{0.75 x},\quad x\in[-1,1],
\end{equation}
which has a critical point of order $\ncp$ at $x = 0$.
Table \ref{tab:ncp0} presents the convergence rates of the WENO-ZA6 scheme with the power parameter $p=1$ for critical points of order $\ncp \le 5$. 
It is observed that the WENO-ZA6 scheme satisfies the sufficient condition in Table \ref{tab:WENO6_condition_ZA_p} for $\ncp=1, 2$, achieving sixth-order accuracy. 
For $\ncp=3,4,5$, where the sufficient condition is no longer satisfied, the WENO-ZA6 scheme exhibits fifth-, fourth-, and fifth-order accuracy, respectively. 
Table \ref{tab:ncp1} summarizes the corresponding results with $p=\frac{3}{2}$ for $\ncp \le 3$ . 
Under this choice, the sufficient condition for sixth-order accuracy in Table \ref{tab:WENO6_condition_ZA_p} is satisfied for all three cases. 
The numerical observations are in complete agreement with the theoretical analysis and the criterion outlined in Table \ref{tab:WENO6_condition_ZA_p}. 
We also show the convergence rates for $\ncp=4,5$ with $p=1$, using the variable sensitivity parameter $\epsilon=\Delta x^{2(\ncp-1)}$ following the aforementioned strategy.
The results clearly demonstrate that this remedy successfully restores sixth-order accuracy even in the presence of higher-order critical points.
\begin{table}[h!]
\centering
\caption{$L_{\infty}$ error and numerical order of accuracy $\order({\dx^m})$ of the WENO-ZA6 scheme ($p=1$) for function \eqref{eq:func_convergence_critical_pt}.}
\label{tab:ncp0}
\renewcommand{\arraystretch}{1.3}
\begin{tabular}{rlcrlcrlcrlcrlc}
\toprule
$\ncp$ & \multicolumn{2}{c}{$1$} & & \multicolumn{2}{c}{$2$} & & \multicolumn{2}{c}{$3$} & & \multicolumn{2}{c}{$4$} & & \multicolumn{2}{c}{$5$} \\
         \cline{2-3}                 \cline{5-6}                 \cline{8-9}                 \cline{11-12}               \cline{14-15}
$1/\dx$  & Error &  $m$          & & Error & $m$             & & Error & $m$             & & Error & $m$             & & Error & $m$ \\
\midrule
20     & 3.0E-09 &     & & 8.4E-06 &     & & 4.2E-05 &     & & 3.5E-05 &     & & 1.6E-05 &     \\
40     & 4.6E-11 & 6.0 & & 1.3E-07 & 6.0 & & 1.8E-06 & {\red 4.6} & & 1.6E-06 & {\red 4.4} & & 4.7E-07 & {\red 5.1} \\
80     & 7.3E-13 & 6.0 & & 1.9E-09 & 6.0 & & 6.4E-08 & {\red 4.8} & & 9.7E-08 & {\red 4.1} & & 1.4E-08 & {\red 5.1} \\
160    & 1.1E-14 & 6.0 & & 3.0E-11 & 6.0 & & 2.1E-09 & {\red 5.0} & & 6.0E-09 & {\red 4.0} & & 4.3E-10 & {\red 5.0} \\
320    & 1.8E-16 & 6.0 & & 4.7E-13 & 6.0 & & 6.4E-11 & {\red 5.0} & & 3.7E-10 & {\red 4.0} & & 1.3E-11 & {\red 5.0} \\
640    & 2.8E-18 & 6.0 & & 7.3E-15 & 6.0 & & 2.0E-12 & {\red 5.0} & & 2.3E-11 & {\red 4.0} & & 4.1E-13 & {\red 5.0} \\
\bottomrule
\end{tabular}
\end{table}

\begin{table}[h!]
\centering
\caption{$L_{\infty}$ error and numerical order of accuracy $\order({\dx^m})$ of the WENO-ZA6 operator ($p = \frac{3}{2}$) for $\ncp = 1,2,3$ and ($p = 1$) for $\ncp = 4,5$ with $\epsilon = \Delta x^{2(\ncp-1)}$.}
\label{tab:ncp1}
\renewcommand{\arraystretch}{1.3}
\begin{tabular}{rlcrlcrlcrlcrlc}
\toprule
$\ncp$ & \multicolumn{2}{c}{$1$} & & \multicolumn{2}{c}{$2$} & & \multicolumn{2}{c}{$3$} & & \multicolumn{2}{c}{$4$} & & \multicolumn{2}{c}{$5$} \\
         \cline{2-3}                 \cline{5-6}                 \cline{8-9}                 \cline{11-12}               \cline{14-15}
$1/\dx$  & Error &  $m$          & & Error & $m$             & & Error & $m$             & & Error & $m$             & & Error & $m$ \\
\midrule
20     & 3.0E-09 &     & & 7.6E-07 &     & & 7.6E-05 &     & & 5.7E-06 &     & & 1.9E-05 &    \\
40     & 4.6E-11 & 6.0 & & 2.9E-09 & 8.0 & & 2.0E-06 & 5.2 & & 1.5E-08 & 8.6 & & 1.6E-07 & 6.9 \\
80     & 7.3E-13 & 6.0 & & 1.2E-11 & 7.9 & & 3.8E-08 & 5.7 & & 2.3E-10 & 6.0 & & 1.5E-09 & 6.7 \\
160    & 1.1E-14 & 6.0 & & 9.0E-14 & 7.1 & & 6.0E-10 & 6.0 & & 3.6E-12 & 6.0 & & 1.7E-11 & 6.5 \\
320    & 1.8E-16 & 6.0 & & 1.4E-15 & 6.0 & & 9.3E-12 & 6.0 & & 5.8E-14 & 6.0 & & 2.3E-13 & 6.0 \\
640    & 2.8E-18 & 6.0 & & 2.2E-17 & 6.0 & & 1.4E-13 & 6.0 & & 8.7E-16 & 6.0 & & 4.2E-15 & 6.0 \\
\bottomrule
\end{tabular}
\end{table}

\subsection{Approximate dispersion relation analysis} 
\label{sec:ADR}

While such a highly symmetric construction is advantageous for achieving high-order accuracy in smooth regions, it may be less favorable for ENO-type shock-capturing and may adversely affect numerical stability near discontinuities.
To elucidate the leading-order nonlinear effects of the proposed shock-capturing scheme, we employ the approximate dispersion relation (ADR) analysis \cite{Pirozzoli}.

\begin{figure}[htb!]
\begin{center}
\mbox{
\makebox[0.4\textwidth][c]{\hskip5pt {\footnotesize Dispersion}}
\makebox[0.4\textwidth][c]{\hskip5pt {\footnotesize Dissipation}}
}
\\[-0.5pt]
\mbox{
\includegraphics[width=0.4\textwidth, trim=10 0 20 10, clip=true]{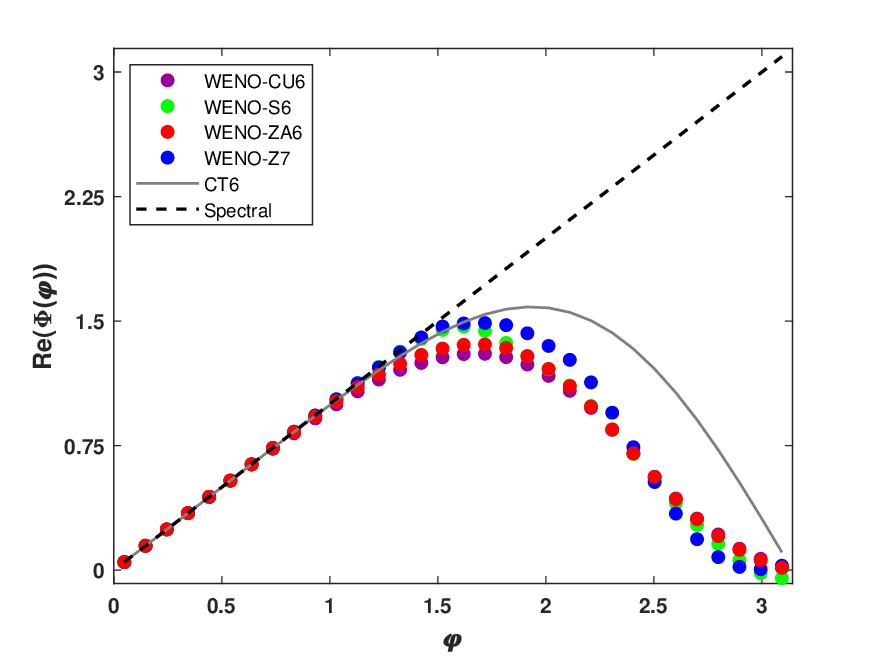}
\includegraphics[width=0.4\textwidth, trim=10 0 20 10, clip=true]{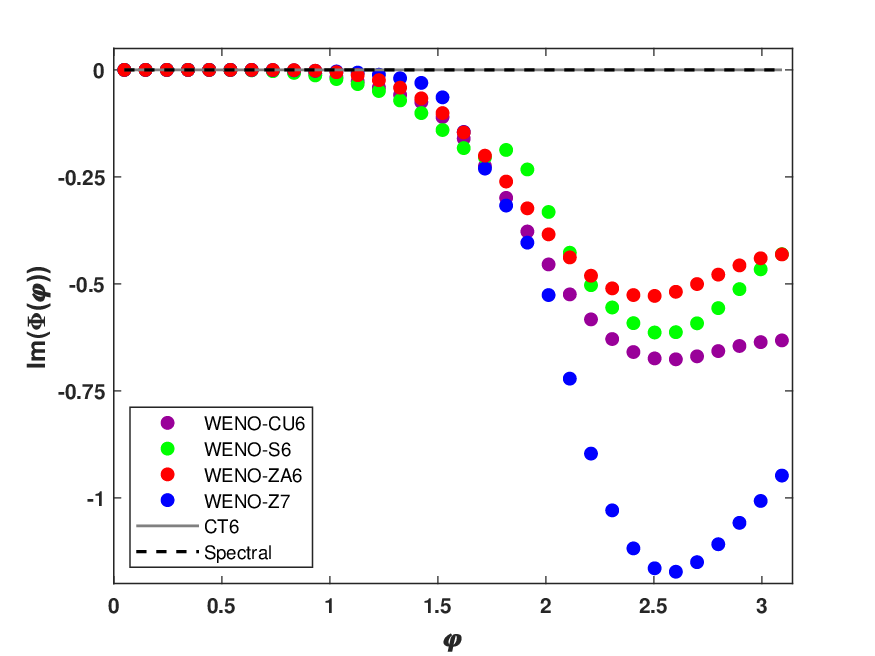}
}
\end{center}
\vskip-20pt
\caption{ADR comparison: $\Real(\Phi(\varphi))$ (phase) versus $\varphi$ (Left) and $\Imag(\Phi(\varphi))$ versus $\varphi$ (Right).
Curves: WENO-CU6 (purple), WENO-S6 (green), WENO-ZA6 (red), WENO-Z7 (blue), CT6 (grey solid), spectral reference (black dashed).}
\label{fig:ADR}
\end{figure}

\begin{figure}[htb!]
\begin{center}
\mbox{
\makebox[0.4\textwidth][c]{\hskip5pt {\footnotesize Dispersion}}
\makebox[0.4\textwidth][c]{\hskip5pt {\footnotesize Dissipation}}
}
\\[-0.5pt]
\mbox{
\includegraphics[width=0.4\textwidth, trim=10 0 20 10, clip=true]{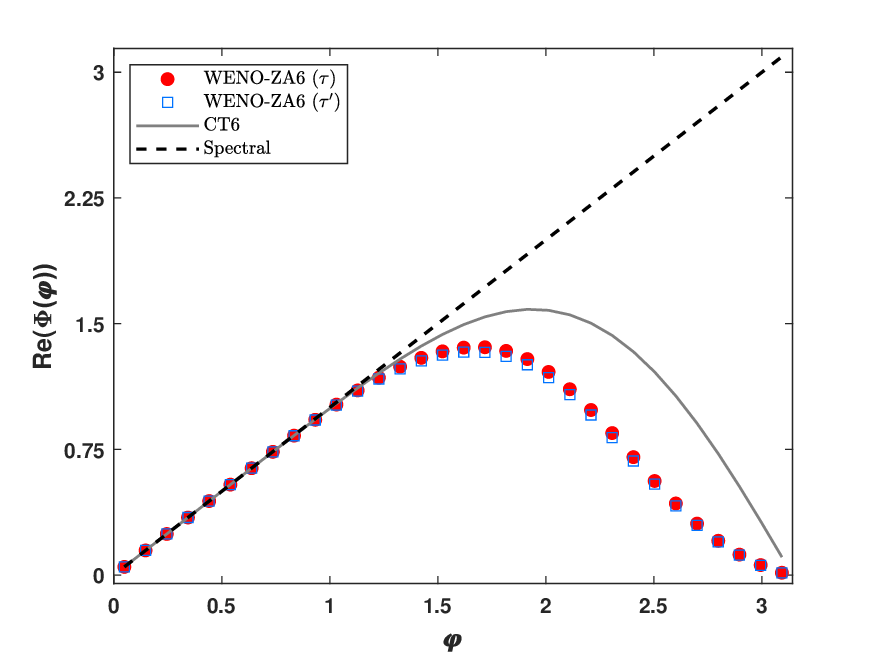}
\includegraphics[width=0.4\textwidth, trim=10 0 20 10, clip=true]{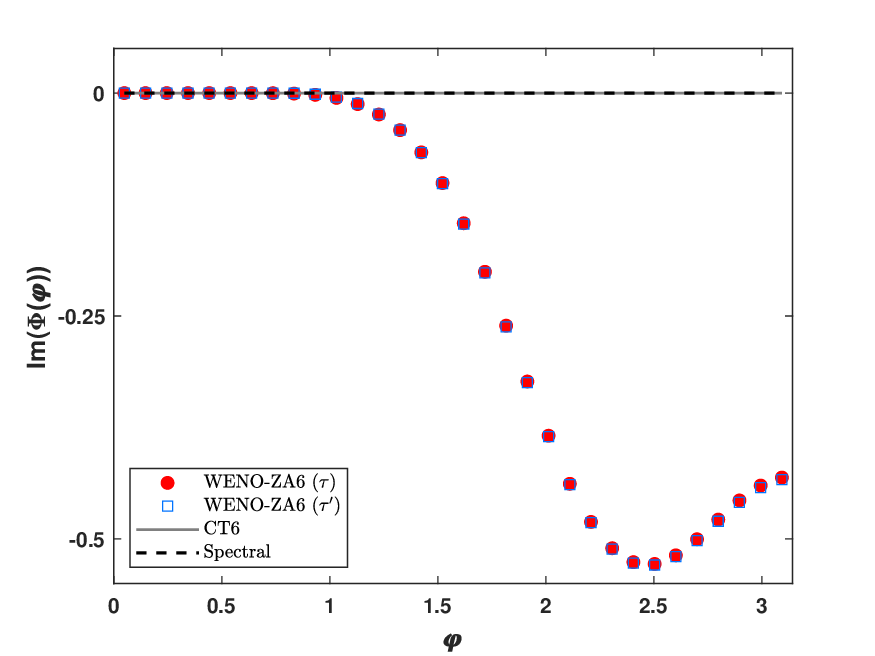}
}\\[0pt]
\mbox{
\includegraphics[width=0.4\textwidth, trim=10 0 20 10, clip=true]{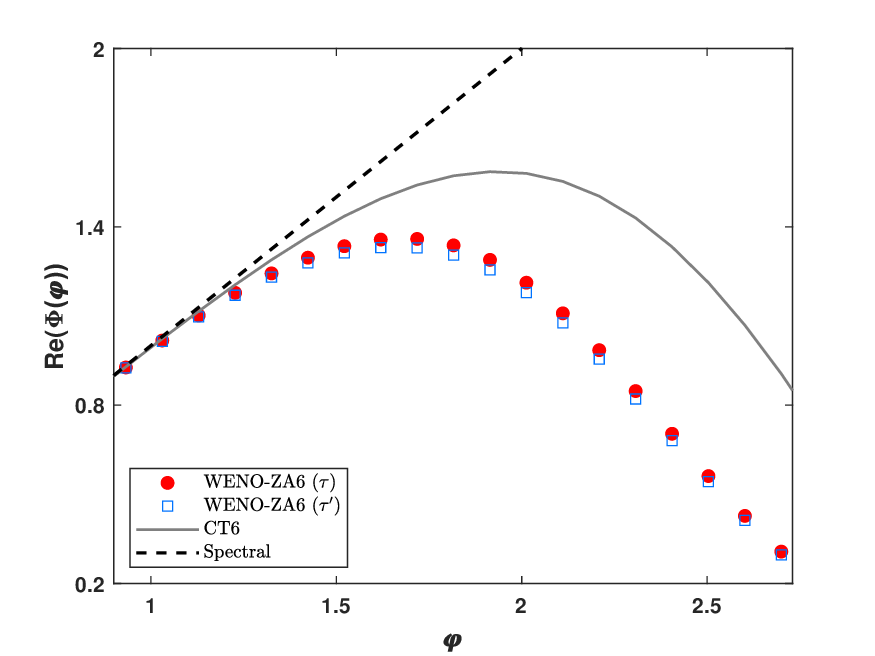}
\includegraphics[width=0.4\textwidth, trim=10 0 20 10, clip=true]{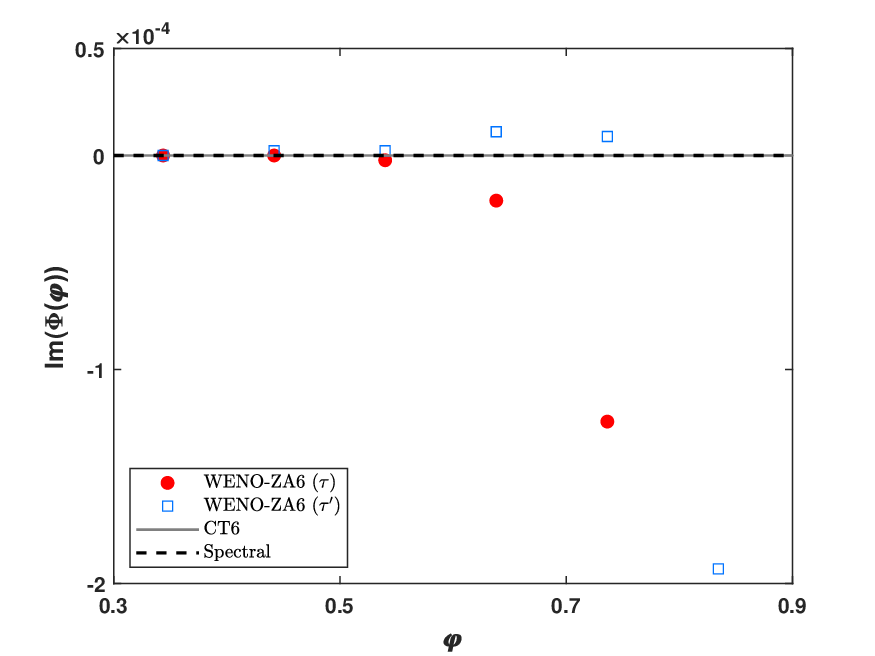}
}
\end{center}
\vskip-20pt
\caption{ADR comparison: $\Real(\Phi(\varphi))$ (phase) versus $\varphi$ (Left) and $\Imag(\Phi(\varphi))$ versus $\varphi$ (Right).
Curves: WENO-ZA6 with $\tau$ \eqref{eq:GSI} (red solid circle), WENO-ZA6 with $\tau'$ \eqref{eq:GSI_alt} (azure hollow square), CT6 (grey solid), spectral reference (black dashed).}
\label{fig:ADR_ZA}
\end{figure}

For each frequency $\varphi$ (or wavenumber $w$ such that $|w \dx|\le \pi$) , we present the real part $\Re\{\Theta(\varphi)\}$ as a measure of dispersion and the imaginary part $\Im\{\Theta(\varphi)\}$ as a measure of dissipation, both plotted against $\varphi$.
The deviation of $\Re\{\Theta(\varphi)\}$ from $\varphi$ quantifies the phase error, whereas the magnitude and structure of $\Im\{\Theta(\varphi)\}$, especially near $\varphi \approx \pi$, indicate the effectiveness of high-wavenumber damping. 
For nonlinear schemes, amplitude dependence in $\Theta(\varphi)$ signals a loss of adaptivity; if the ADR collapses to that of a linear scheme, then the nonlinear weights are no longer effective.
In Fig. \ref{fig:ADR}, we compare the exact spatial differentiation (spectral), the sixth-order linear central (CT6), WENO-CU6 \cite{HuWangAdams}, WENO-S6 \cite{Huang}, WENO-Z7 \cite{DonBorges}, and the proposed WENO-ZA6 schemes. 
For WENO-ZA6, we also consider both global smoothness indicators $\tau$ \eqref{eq:GSI} and $\tau'$ \eqref{eq:GSI_alt}.
Fig. \ref{fig:ADR_ZA} shows that WENO-ZA6 with $\tau$ exhibits dispersion characteristics that are closer to those of CT6 than the version based on $\tau'$, indicating a more accurate phase propagation in WENO-ZA6 with $\tau$ for medium-to-high wavenumbers. 
By contrast, WENO-ZA6 with $\tau'$ introduces a region of negative numerical dissipation around $\varphi = 0.7$, which may lead to spurious energy amplification and reduced robustness. 
These results suggest that $\tau$ should be chosen in the WENO-ZA6 framework, as it plays a critical role in balancing dispersion accuracy and numerical stability.

\section{Numerical experiments}
\label{sec:numerics}

We present one- and two-dimensional numerical results to demonstrate the performance of the proposed WENO-ZA6 scheme.
For the nonlinear component $\Gamma_k$ in \eqref{eq:weight_nonlinear}, $p=1$ is used as the sufficient condition \eqref{eq:WENO6_condition} is satisfied for $\ncp=1, 2$. 
We compare the numerical solutions from WENO-ZA6 with those of WENO-CU6 \cite{HuWangAdams}, WENO-S6 \cite{Huang}, and WENO-Z7 \cite{DonBorges}.
The sensitivity parameter is chosen as $\eps = 10^{-40}$ for WENO-ZA6, while the original settings of $\eps$ are retained for the other three WENO schemes.
We integrate in time using the explicit third-order total variation diminishing Runge-Kutta method \cite{ShuOsherI} with $\cfl = 0.45$.

\subsection{One-dimensional scalar equation} \label{sec:scalar1d}
\begin{figure}[htb!]
\begin{center}
\mbox{
\includegraphics[width=0.33\textwidth, trim=5 5 5 10, clip=true]{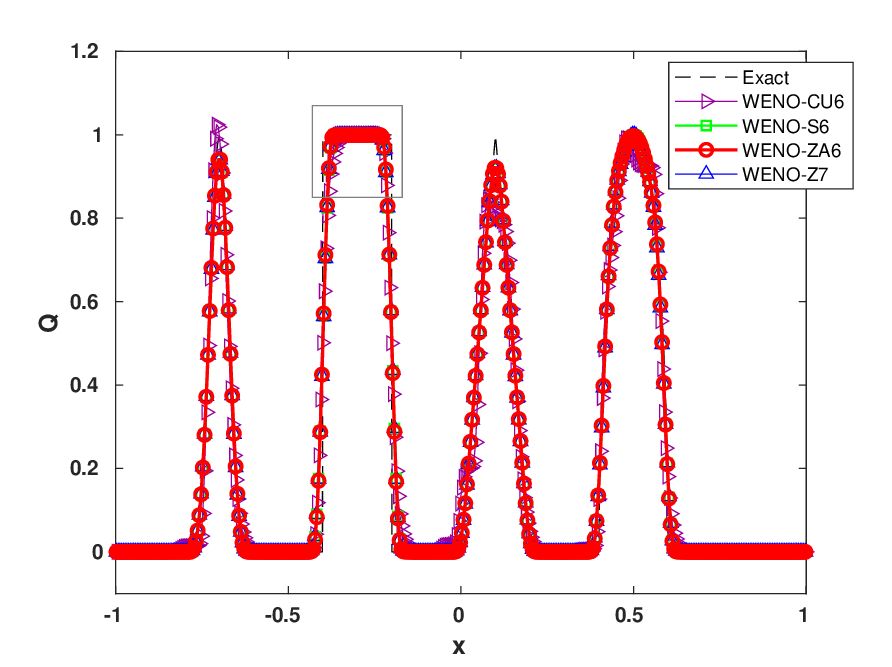}
\includegraphics[width=0.33\textwidth, trim=5 5 5 10, clip=true]{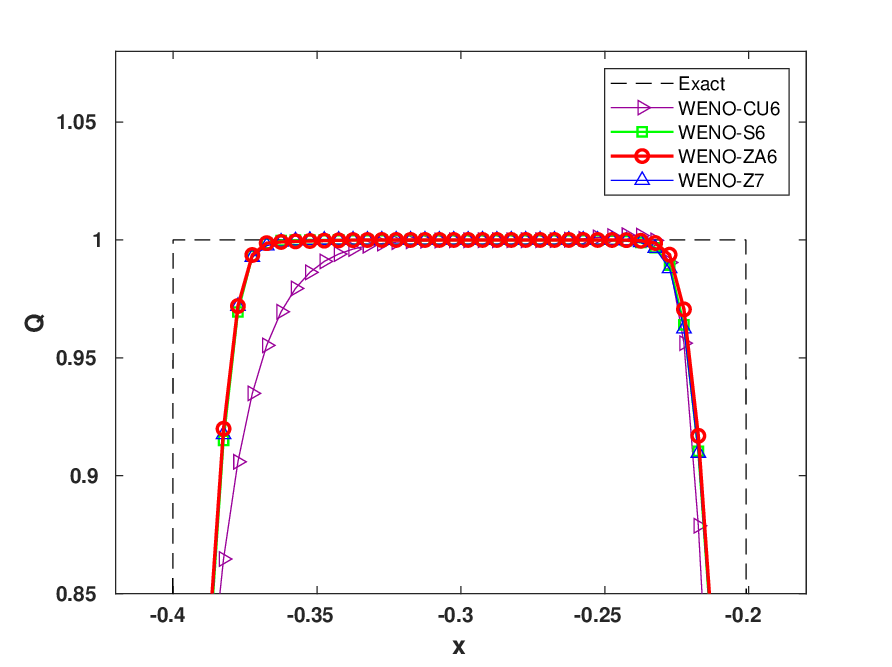}
\includegraphics[width=0.33\textwidth, trim=5 5 5 10, clip=true]{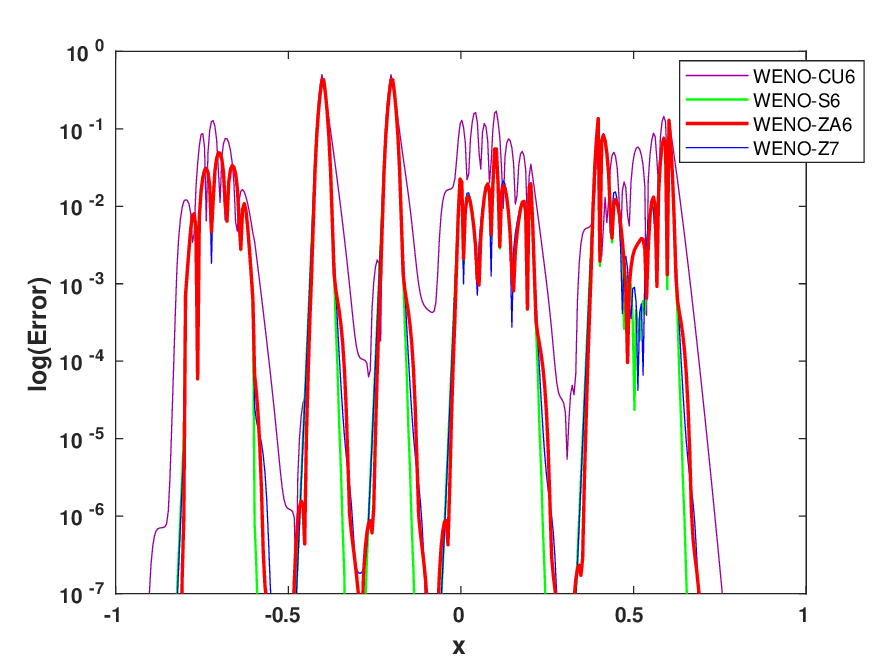}
}
\end{center}
\vskip-10pt
\caption{Example \ref{ex:multiwave}, (Left) numerical solution, (Middle) close-up view of the square wave, and (Right) log-scale absolute pointwise error of multi-waves advection problem computed by the WENO-CU6 (purple), WENO-S6 (green), WENO-ZA6 (red), and WENO-Z7 (blue) schemes with $N = 400$ at $T = 20$.}
\label{fig:multiwave}
\end{figure}

\begin{example}\label{ex:multiwave}\rm 
In this test, we apply the WENO schemes to the multi-wave advection equation $Q_t+ Q_x = 0$ with the initial condition consisting of a Gaussian, triangle, square wave, and semi-ellipse, given by 
\begin{align*}\label{eq:multiwave}
& Q(x,0) = \left\{ 
            \begin{array}{lc} 
             \frac{1}{6} \left[ G(x, \beta, z-\delta) + 4 G(x, \beta, z) + G(x, \beta, z+\delta) \right],   & -0.8 \leq x \leq -0.6, \\
             1,                                                                                             & -0.4 \leq x \leq -0.2, \\
             1 - \left| 10(x-0.1) \right|,                                                                  & 0 \leq x \leq 0.2, \\
             \frac{1}{6} \left[ F(x, \alpha, y-\delta) + 4 F(x, \alpha, y) + F(x, \alpha, y+\delta)\right], & 0.4 \leqslant x \leqslant 0.6, \\
             0, & \text{otherwise},
            \end{array} 
           \right. \\
& G(x, \beta, z) = \e^{-\beta(x-z)^2}, \; F(x, \alpha, y) = \sqrt{\max \left( 1-\alpha^2 (x-y)^2, 0 \right)}, \\
& \delta = 0.005, \, \beta = \frac{\ln 2}{36 \delta^2}, \, z = -0.7, \, \alpha = 10, \, y = 0.5.
\end{align*}
The final time is $T=20$ on the interval $[-1,1]$ with periodic boundary conditions. 
The number of mesh cells is $N = 400$. 

In Fig.~\ref{fig:multiwave}, the numerical solutions and their zoomed figure around the square wave are shown. Absolute pointwise errors computed by the WENO schemes are also presented. These results illustrate that WENO-S6, ZA6, and Z7 outperform WENO-CU6. The latter (purple line) exhibits stronger numerical dissipation around the square wave and suffers larger spurious oscillations away from the discontinuities in the other shapes.
The absolute pointwise errors decrease monotonically and exponentially away from the singularities (discontinuities and sharp corners) for all WENO schemes. 
This behavior captures the discontinuities in the ENO manner and agrees well with the exact solution, particularly for the square wave. Furthermore, in the smooth regions of the solution, WENO-ZA6 is competitive with WENO-S6 and WENO-Z7. It performs substantially better than WENO-CU6, as indicated in the pointwise error plot.
\end{example}

\begin{example}[Scalar conservation laws with convex and nonconvex fluxes]
\label{ex:scalar_fluxes}\rm
We consider the one-dimensional scalar hyperbolic conservation law \eqref{eq:hyperbolic_1D_scalar} with different choices of the flux function $f(Q)$. 
All problems differ only in the types of flux function (convex or nonconvex), while the governing equation and numerical framework remain the same. 
The detailed setup for each case is summarized in Table~\ref{tab:scalar_fluxes}. 

\begin{table}[htb!]
\begin{center}
\caption{Problem setup for scalar conservation laws with different choices of the flux function.}
\label{tab:scalar_fluxes}
\begin{tabular}{l l l l l}
\hline
& Flux type & $f(Q)$ & $Q(x,0)$ & $(N,\,T)$ \\
\hline

Burgers 
& Convex 
& $\dfrac{Q^2}{2}$ 
& $\dfrac{1}{2} + \sin(2\pi x)$ 
& $(80,\, 0.2)$ \\ [6pt]

Buckley--Leverett 
& Nonconvex 
& $\dfrac{Q^2}{Q^2 + 0.25(1-Q)^2}$ 
& $\begin{cases}
1, & -0.5 \le x \le 0, \\
0, & \text{otherwise}
\end{cases}$
& $(200,\, 0.3)$ \\

\hline
\end{tabular}
\end{center}
\end{table}

\begin{figure}[htb!]
\begin{center}
\mbox{
\makebox[0.33\textwidth][c]{Burgers equation}
\makebox[0.33\textwidth][c]{Pointwise error}
\makebox[0.33\textwidth][c]{Buckley-Leverett problem}
}
\mbox{
\includegraphics[width=0.33\textwidth, trim=10 5 30 10, clip=true]{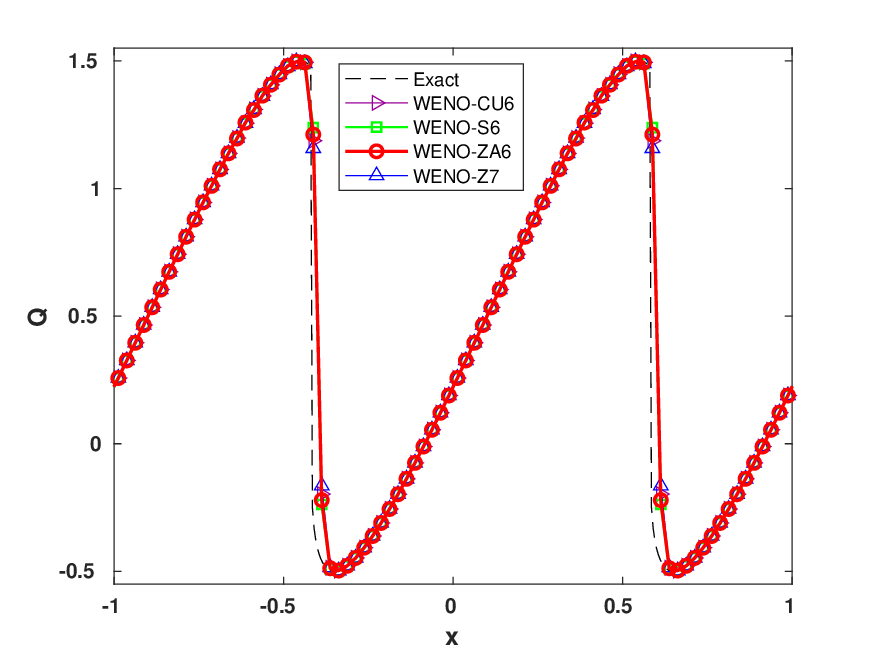}
\includegraphics[width=0.33\textwidth, trim=10 5 30 10, clip=true]{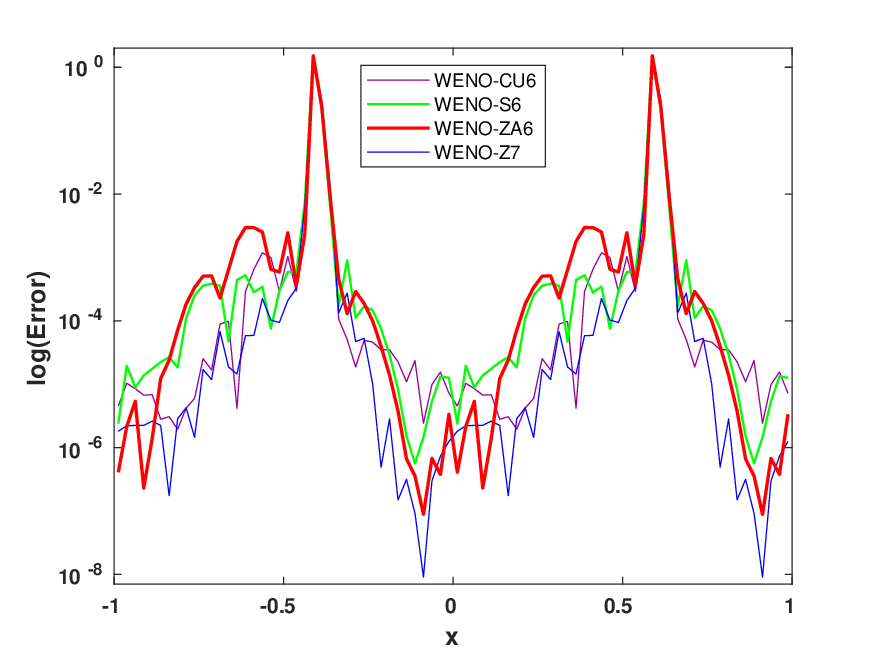}
\includegraphics[width=0.33\textwidth, trim=10 5 30 10, clip=true]{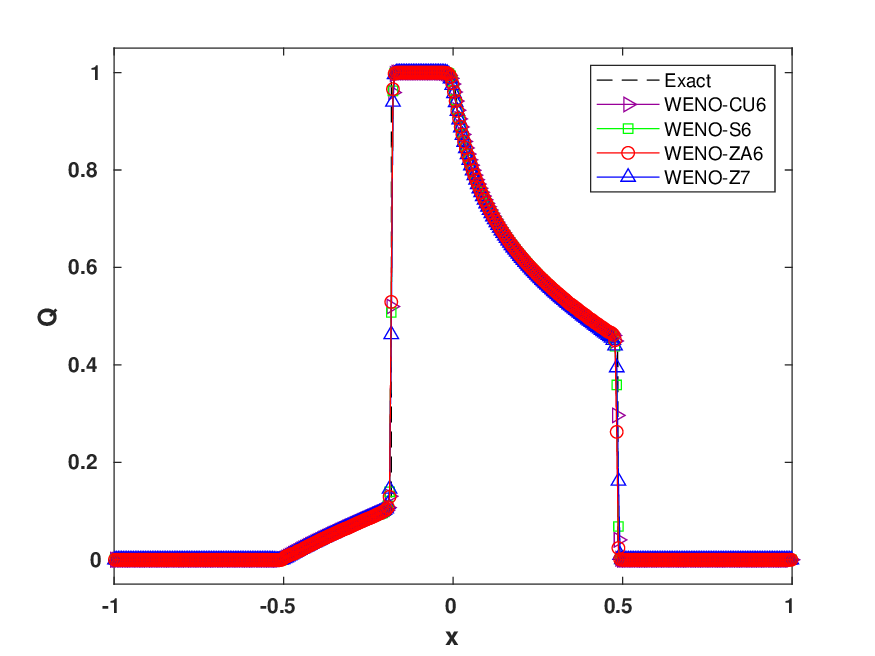}
}
\end{center}
\vskip-10pt
\caption{Example \ref{ex:scalar_fluxes}: (Left) the numerical solution and (Middle) absolute pointwise error of Burgers equation at $T=0.2$, and (Right) the numerical solution of the Buckley-Leverett problem at $T=0.3$ computed by the WENO-CU6 (purple), WENO-S6 (green), WENO-ZA6 (red), and WENO-Z7 (blue) schemes.
}
\label{fig:burgers}
\end{figure}
These test cases assess the robustness and accuracy of the WENO schemes under distinct wave structures induced by the flux function.
For the convex quadratic flux of the Burgers equation, the moving solution remains smooth up to the critical time $\tfrac{1}{2\pi} \approx 0.1592$, where the shock formation begins. 
As shown in the right subfigure of Fig. \ref{fig:burgers}, all schemes produce highly accurate approximations at the final time $T = 0.2$. 
In the middle subfigure, the pointwise errors on a logarithmic scale are uniformly small across the smooth regions of the solution and far from the developing shock, and are $O(1)$ near the neighborhood of a steep gradient. The error decays exponentially and monotonically, showing a sharp, essentially non-oscillatory capture of the forming shock.
The S-shaped flux of the Buckley-Leverett problem generates two pairs of composite wave structures with attached shock and rarefaction waves.
As illustrated in the left subfigure, all schemes successfully resolve the complex wave interactions in a stable and ENO manner. 
Overall, the results demonstrate that the proposed WENO schemes maintain high accuracy in smooth regions and robust shock-capturing capability across convex and nonconvex flux types.
\end{example}

\subsection{One-dimensional Euler equations} \label{sec:Euler1d}
Next, we solve the one-dimensional Euler equations of gas dynamics,
\begin{equation}\label{eq:Euler1D}
 \frac{\partial}{\partial{t}}
\begin{pmatrix}
\rho\\
\rho u\\
E
\end{pmatrix}
 +\frac{\partial}{\partial{x}}
\begin{pmatrix}
\rho u\\
\rho u^2+p\\
(E+p)u
\end{pmatrix} = 0 ,
\end{equation}
where $\rho$, $u$, $E$, and $p$ denote the density, velocity, total energy, and pressure, respectively. 
The pressure is determined by the equation of state
$$ p = (\gamma - 1)\left(E-\half \rho u^2\right). $$
The ratio of specific heats is taken to be $\gamma=1.4$ unless stated otherwise.

We compare the performance of the WENO schemes on several classical shock-tube benchmark problems with Riemann initial data, including the Sod, Lax, and 123 problems, the shock-entropy interaction problem, the shock-density interaction problem, and the two-blast-wave interaction problem.
The characteristic-wise Lax-Friedrichs flux splitting \cite{ShuSpringer} is used for the Euler equations in this subsection and the next.

\begin{figure}[htb!]
\begin{center}
\mbox{
\makebox[0.33\textwidth][c]{\hskip15pt Sod}
\makebox[0.33\textwidth][c]{\hskip15pt Lax}
\makebox[0.33\textwidth][c]{\hskip15pt 123}
}
\mbox{
\includegraphics[width=0.33\textwidth, trim=12 5 39 10, clip=true]{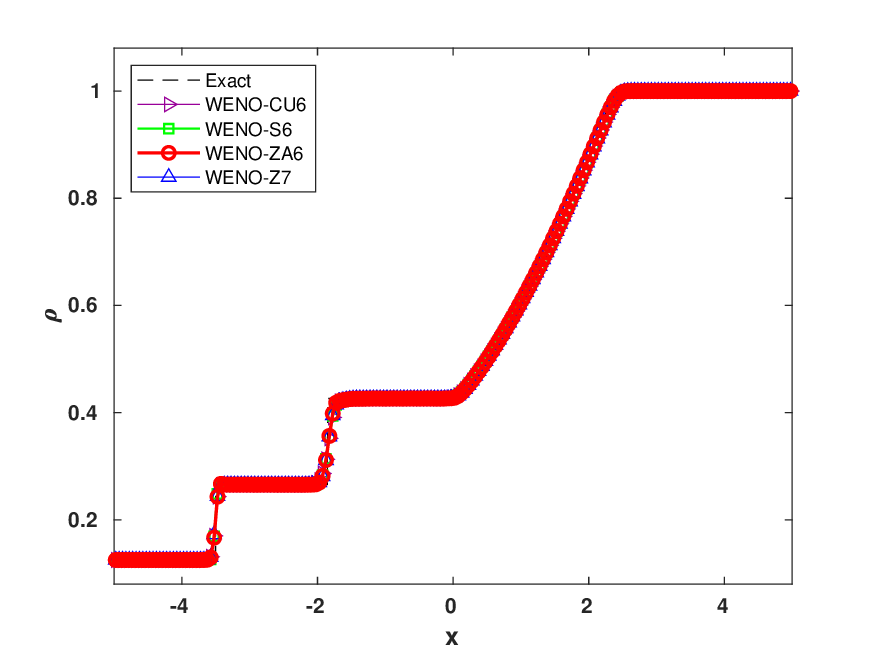}
\includegraphics[width=0.33\textwidth, trim=12 5 39 10, clip=true]{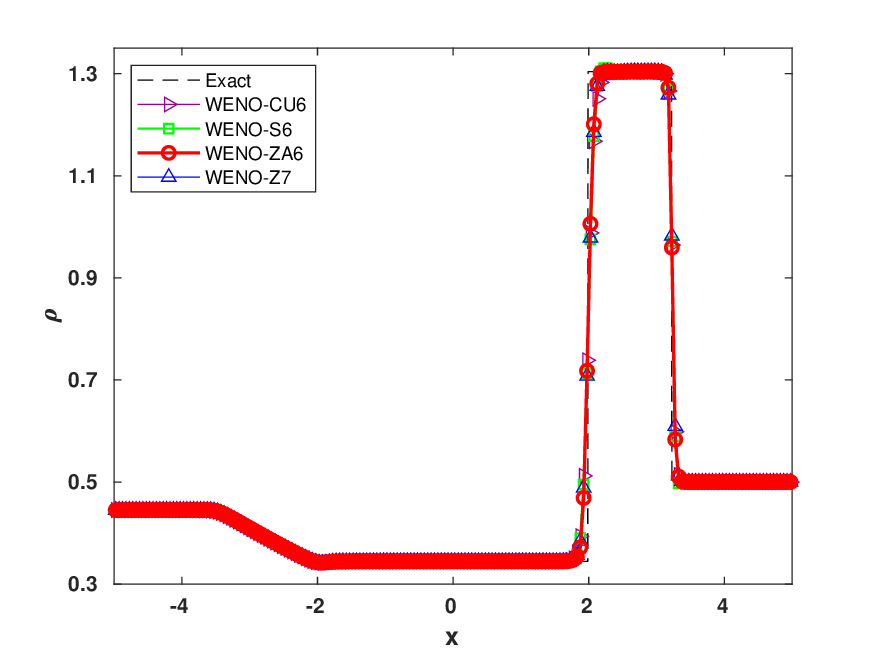}
\includegraphics[width=0.33\textwidth, trim=12 5 39 10, clip=true]{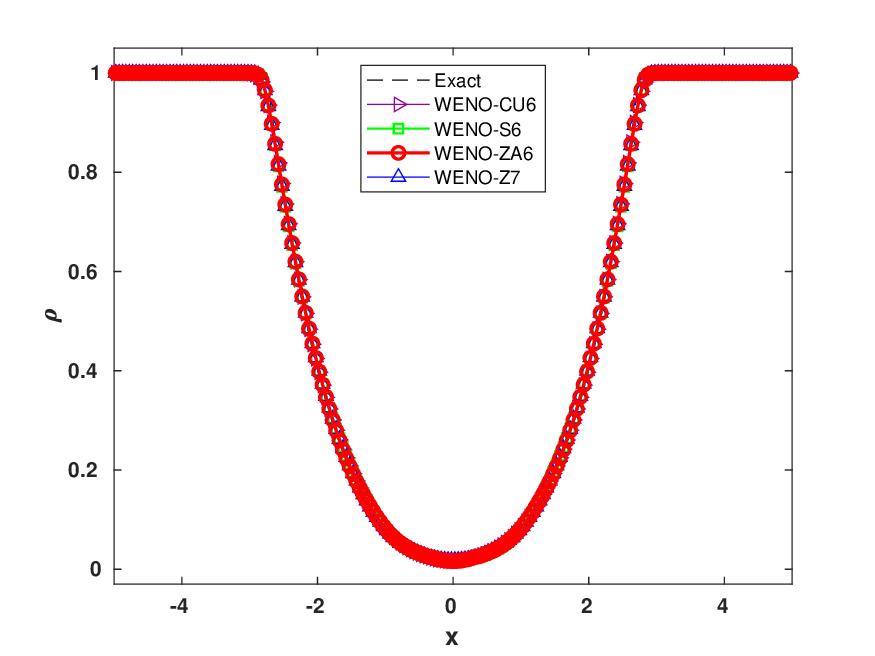}
}
\mbox{
\includegraphics[width=0.33\textwidth, trim=12 5 39 10, clip=true]{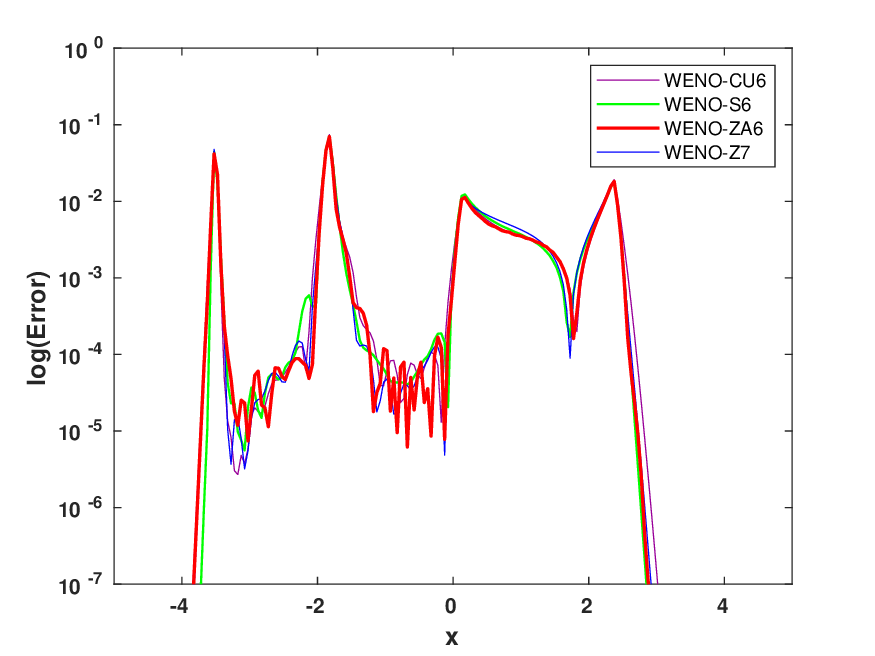}
\includegraphics[width=0.33\textwidth, trim=12 5 39 10, clip=true]{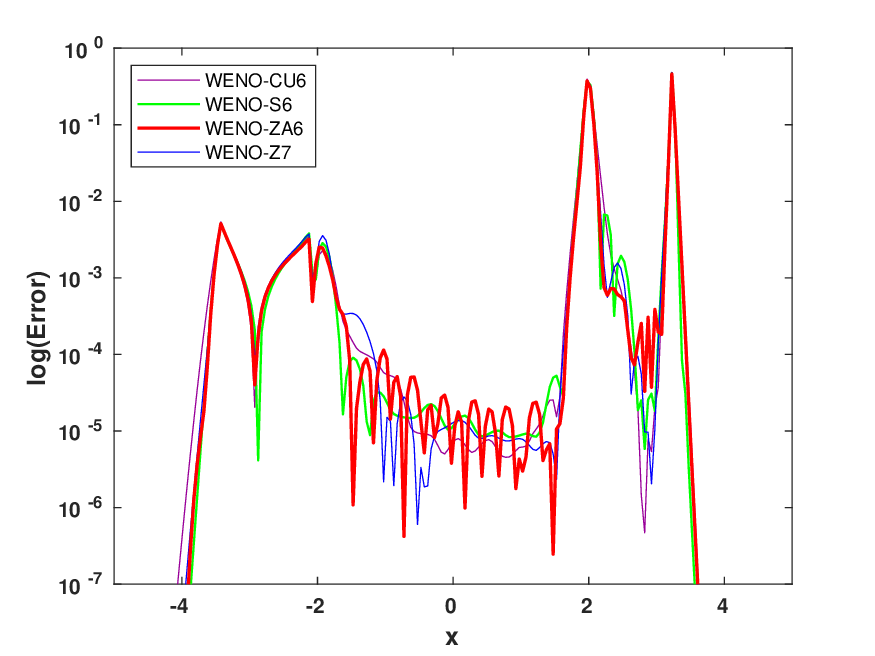}
\includegraphics[width=0.33\textwidth, trim=12 5 39 10, clip=true]{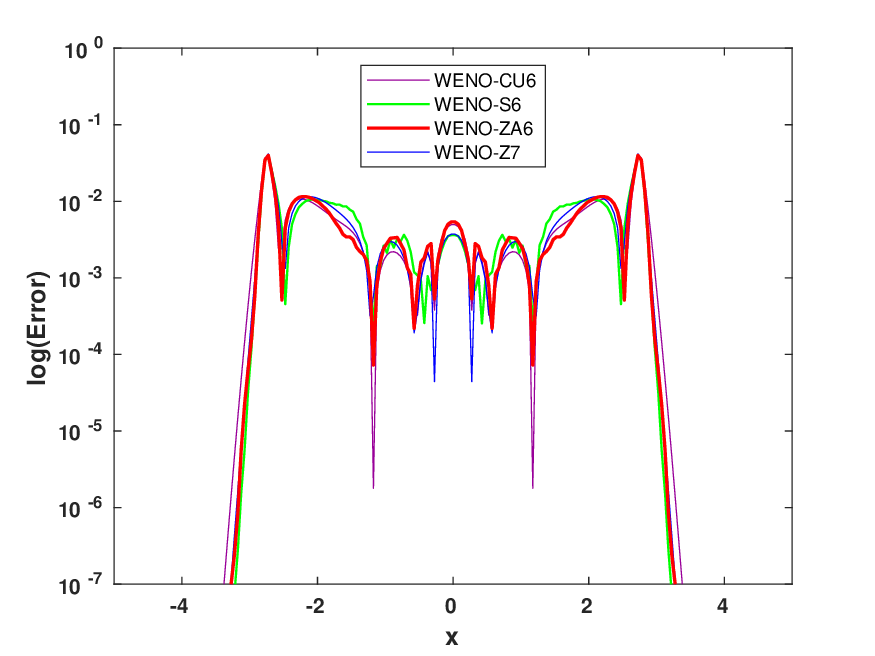}
}
\end{center}
\vskip-10pt
\caption{Example \ref{ex:Riemann}, (Top) density profiles and (Bottom) their absolute pointwise error of (Left) Sod problem $(T=2)$, (Middle) Lax problem $(T=1.3)$ and (Right) 123 problem ($T=1)$ by the WENO-CU6 (purple), WENO-S6 (green), WENO-ZA6 (red), and WENO-Z7 (blue) schemes with $N = 200$.}
\label{fig:Riemann}
\end{figure}

\begin{example}\label{ex:Riemann}\rm 
To demonstrate the ENO shock-capturing capabilities of the proposed WENO-ZA6 scheme, we present three classical shock-tube benchmark problems with high gradients, discontinuities, and small weak smooth structures: the Sod, Lax, and 123 problems. 
The numerical experiments are conducted using $N=200$ mesh cells for the computational domain $[-5, \, 5]$. 
The corresponding initial states and final times are listed in Table \ref{tab:Euler1D_Riemann}, and the shock locations are set as $x=0$ initially.

\begin{table}[htb!]
\centering
\caption{Problem setup for one-dimensional Riemann problems.}
\label{tab:Euler1D_Riemann}
\begin{tabular}{l c c c c c c c c}
\hline
& $\rho_L$ & $u_L$ & $p_L$ &
& $\rho_R$ & $u_R$ & $p_R$
& T \\
\hline

Sod
& 0.125 & 0 & 0.1 &
& 1     & 0 & 1
& 2 \\

Lax
& 0.445 & 0.698 & 3.528 &
& 0.5   & 0     & 0.571 
& 1.3 \\

123
& 1 & -2 & 0.4 &
& 1 &  2 & 0.4
& 1 \\
\hline
\end{tabular}
\end{table}

Figure \ref{fig:Riemann} displays the density $\rho$ and its absolute pointwise error of the Sod, Lax, and 123 problems. 
Numerical results reveal that WENO schemes yield exceptional solutions, efficiently resolving large-scale structures and achieving sharp capture of the contact discontinuity and shock wave, in an ENO manner by virtue of monotone and faster exponential decay of the pointwise error around the singularities.  
\end{example}

\begin{figure}[htb!]
\begin{center}
\mbox{
\includegraphics[width=0.4\textwidth, trim=10 5 30 10, clip=true]{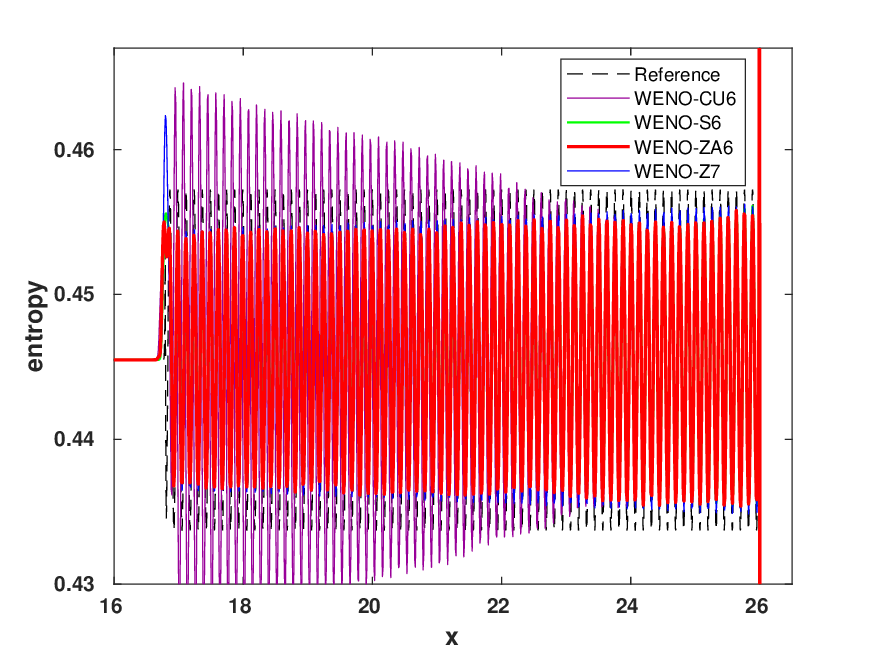}~~
\includegraphics[width=0.4\textwidth, trim=10 5 30 10, clip=true]{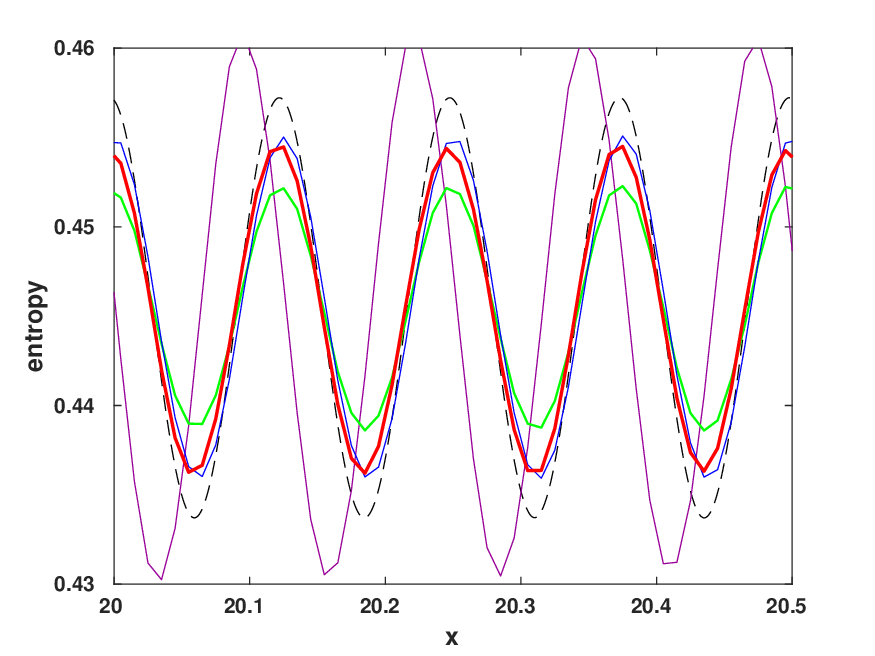}
}
\end{center}
\vskip-10pt
\caption{Example \ref{ex:entropy}, (Left) numerical entropy of shock-entropy wave interaction problem and (Right) close-up view of the entropy by the WENO-CU6 (purple), WENO-S6 (green), WENO-ZA6 (red), and WENO-Z7 (blue) schemes with $N = 4000$ at $T = 10$.}
\label{fig:entropy}
\end{figure}

\begin{example}\label{ex:entropy}\rm
Consider a long-time right-moving Mach 3 shock interacting with a small amplitude sinusoidal perturbation of the entropy in the pre-shock region. 
The initial condition is
$$
\left(\rho,~u,~P\right) = 
\left\{
\begin{array}{lrrrrr}
(& \frac{27}{7},                &\frac{4\sqrt{35}}{9}, & \frac{31}{3} &), & x \leq -9.5, \\
(& \exp(- 0.01\sin(13(x-9.5))), & 0,                   & 1            &), & x > -9.5.
\end{array}
\right.
$$
We take a uniform mesh with $N=4000$ cells over the computational domain $[-10, \, 30]$ and the final time $T=10$ to enable continuous generation of high-frequency waves behind the main shock, thereby testing the performance of the WENO schemes in a long-time simulation. 
Figure \ref{fig:entropy} shows the entropy of the Mach 3 shock-entropy wave interaction computed by the WENO schemes. 
As expected, high-frequency entropy waves computed by WENO-ZA6 and WENO-Z7 are comparable to the reference solution, which is generated by WENO-Z7 with $N = 40000$. 
It is also evident that the corresponding waves computed by WENO-CU6 exhibit noticeable phase errors and antidiffusion.
These results indicate that, on this mesh, WENO-ZA6 and WENO-Z7 achieve higher resolution and lower dissipation than WENO-S6 in long-time simulations.
\end{example}

\begin{figure}[htb!]
\begin{center}
\mbox{
\includegraphics[width=0.33\textwidth, trim=15 5 35 10, clip=true]{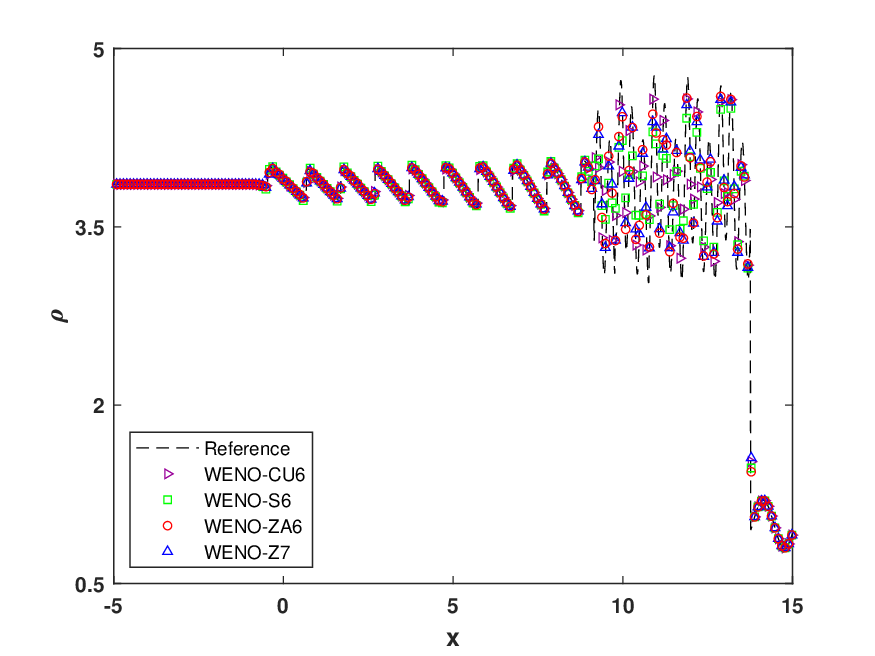}
\includegraphics[width=0.33\textwidth, trim=15 5 35 10, clip=true]{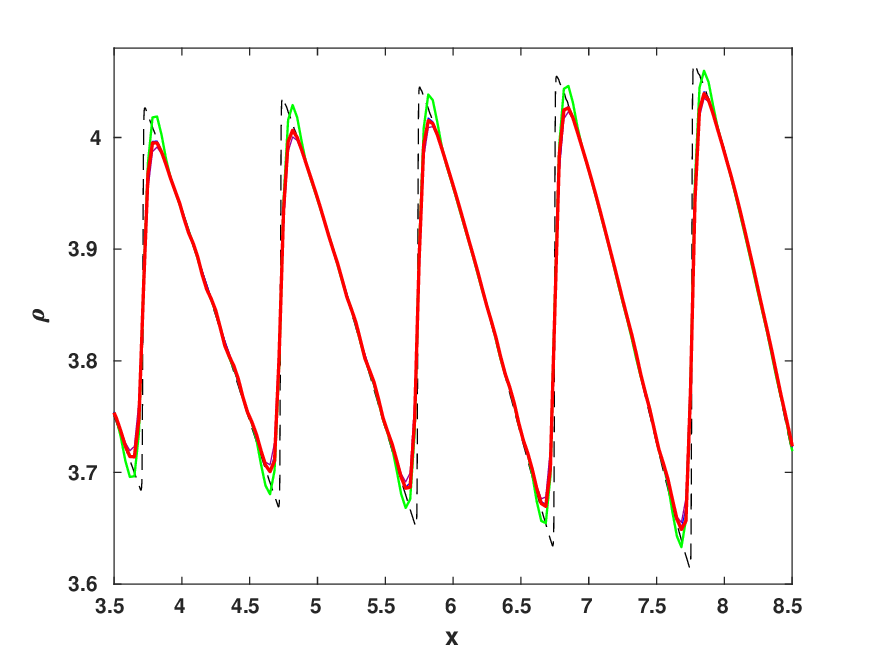}
\includegraphics[width=0.33\textwidth, trim=15 5 35 10, clip=true]{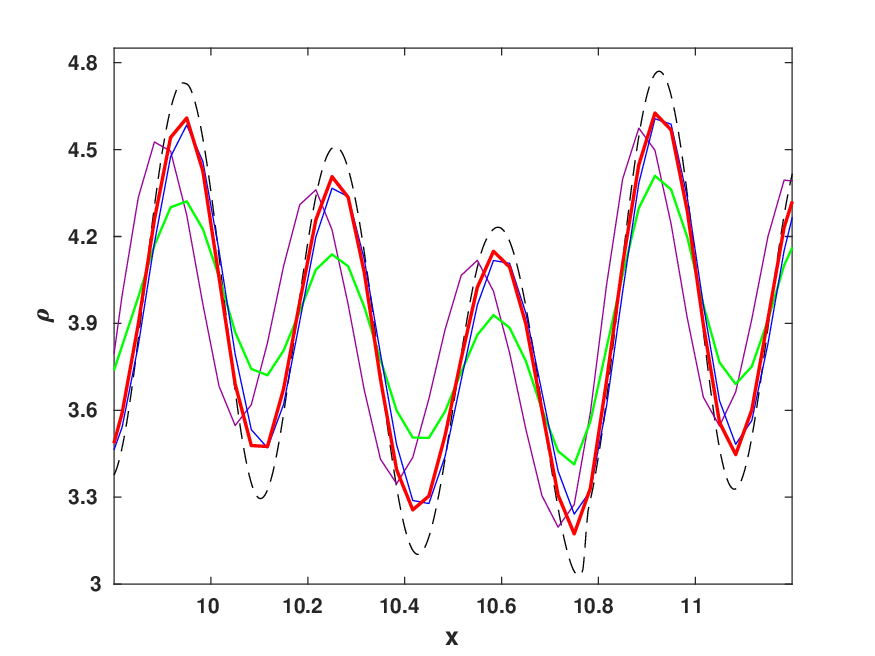}
}
\mbox{
\includegraphics[width=0.245\textwidth, trim=15 5 30 10, clip=true]{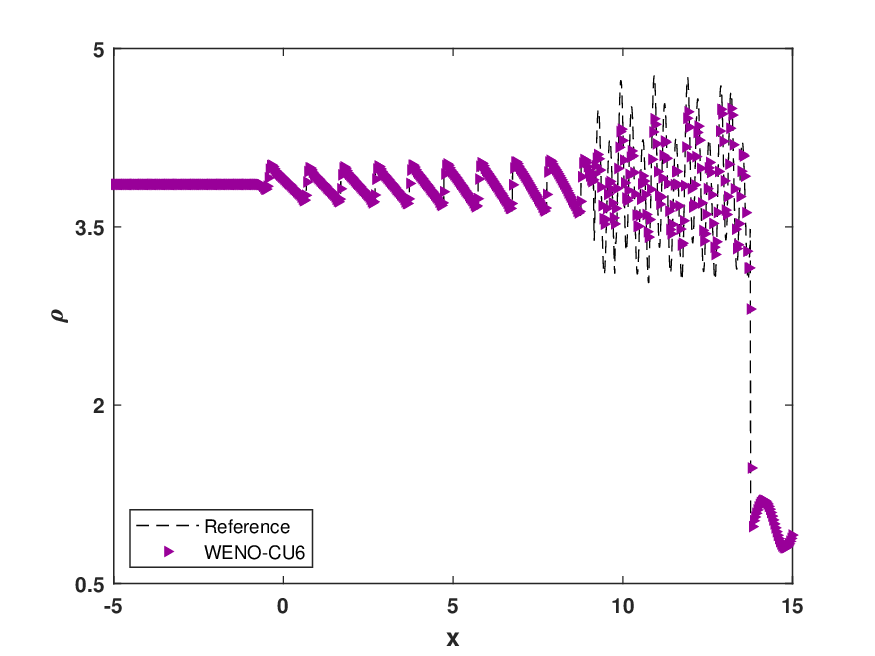}
\includegraphics[width=0.245\textwidth, trim=15 5 30 10, clip=true]{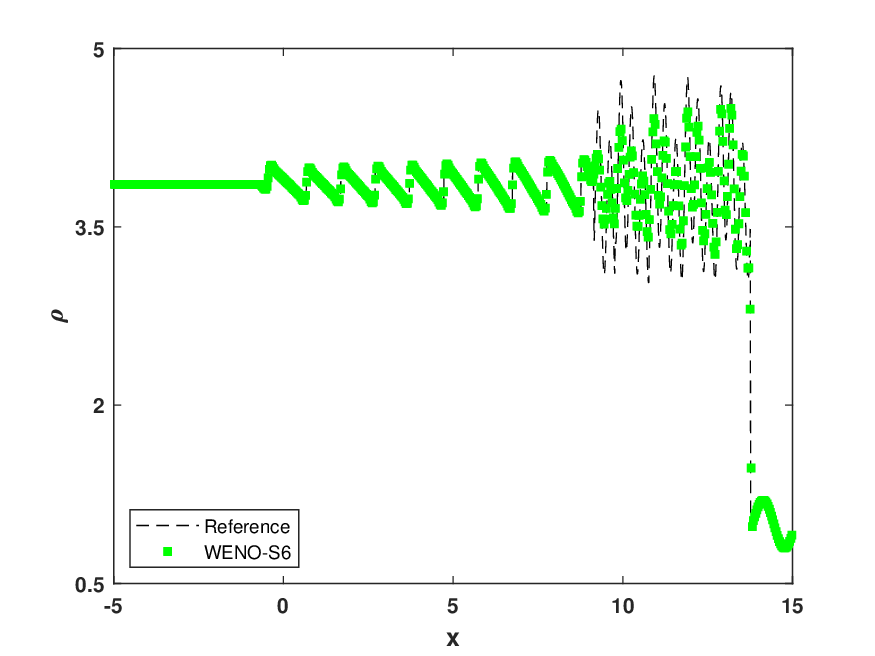}
\includegraphics[width=0.245\textwidth, trim=15 5 30 10, clip=true]{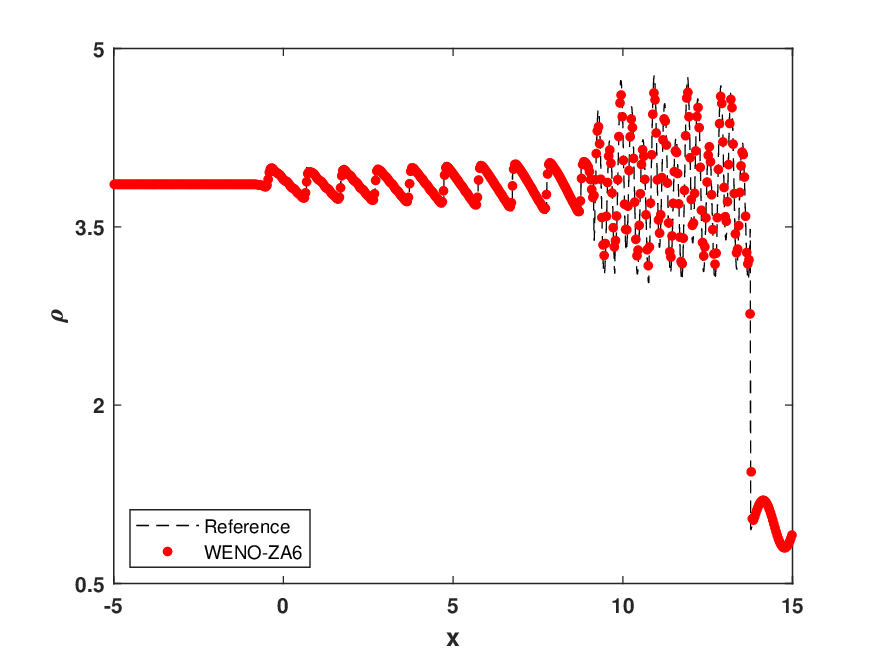}
\includegraphics[width=0.245\textwidth, trim=15 5 30 10, clip=true]{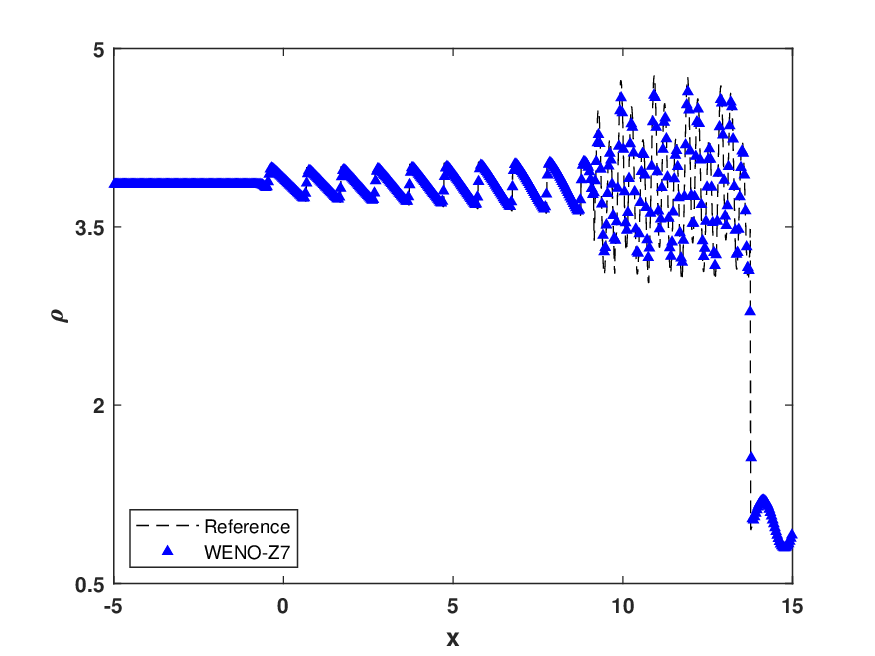}
}
\end{center}
\vskip-10pt
\caption{Example \ref{ex:density}, (Left) numerical density of shock-density wave interaction problem and (Middle, Right) close-up view of the density by the WENO-CU6 (purple), WENO-S6 (green), WENO-ZA6 (red), and WENO-Z7 (blue) schemes with $N = 600$ at $T=5$.
}
\label{fig:density}
\end{figure}

\begin{example}\label{ex:density}\rm 
The extended Mach 3 shock-density wave interaction problem, which solves smooth high-frequency waves behind a right-going main shock, has the initial condition
$$
\left(\rho,~u,~P\right) = 
\left\{
\begin{array}{lrrrrl}
(& \frac{27}{7},  &\frac{4\sqrt{35}}{9}, & \frac{31}{3} &), & x \leq -4, \\
(& 1+0.2\sin(5x), & 0,                   & 1            &), &  x > -4.
\end{array}
\right.
$$
We divide the computational domain $[-5, 15]$ into $N=600$ uniform cells to observe the continuous creation and long-time evolution of shocklets and high-frequency waves in space and time. 
Figure~\ref{fig:density} shows the density profiles computed by the considered WENO schemes at \(T = 5\). 
The reference solution is obtained using WENO-Z7 on a fine grid with \(N = 6000\) cells. 
All schemes accurately capture both the main and emerging shocklets. 
Moreover, WENO-ZA6 performs comparably to WENO-Z7 and clearly outperforms WENO-CU6 and WENO-S6 in resolving the high-frequency acoustic waves.

\end{example}

\begin{figure}[htb!]
\begin{center}
\mbox{
\includegraphics[width=0.4\textwidth, trim=30 4 30 10, clip=true]{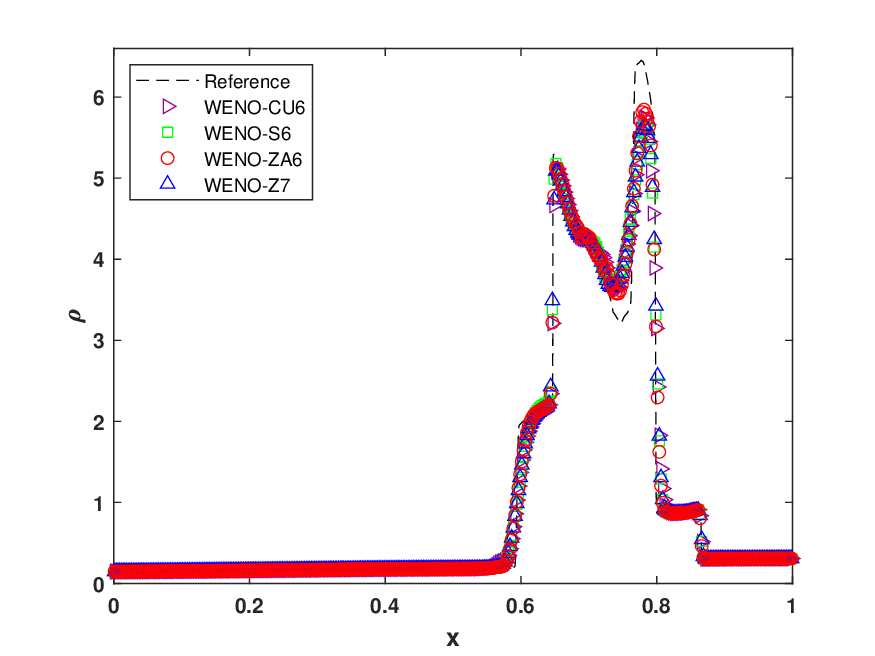}
\includegraphics[width=0.4\textwidth, trim=30 4 30 10, clip=true]{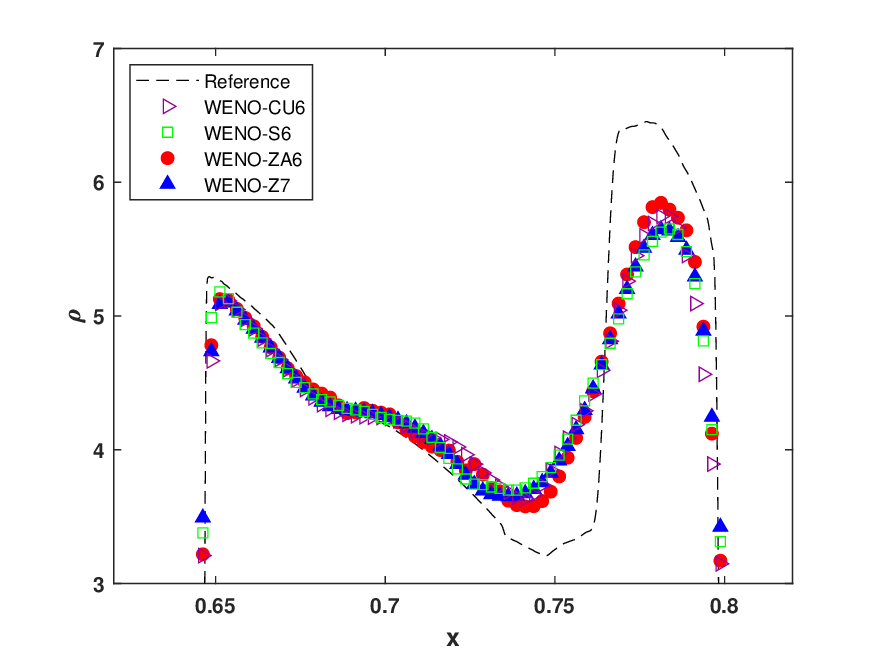}
}
\mbox{
\includegraphics[width=0.245\textwidth, trim=30 5 30 10, clip=true]{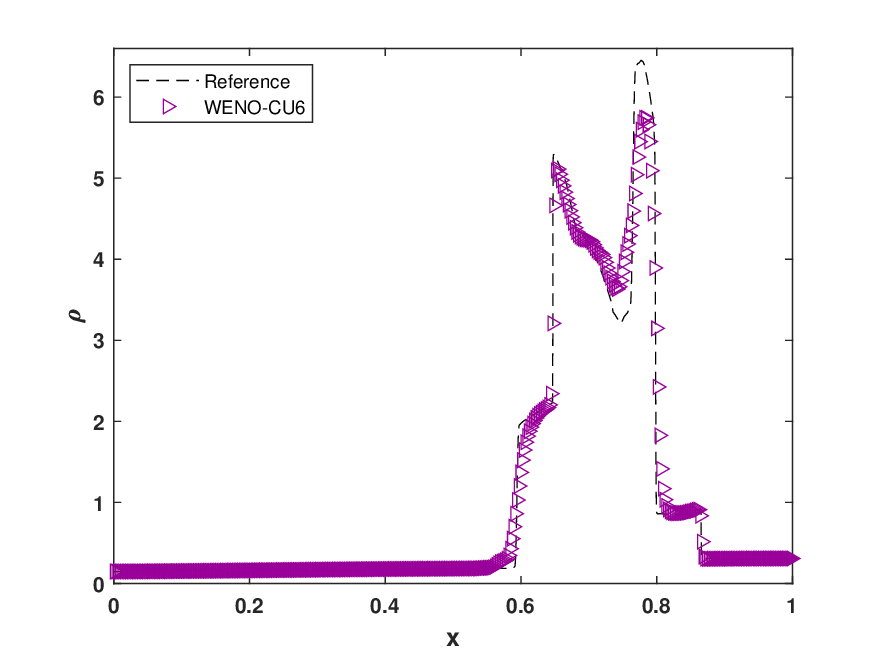}
\includegraphics[width=0.245\textwidth, trim=30 5 30 10, clip=true]{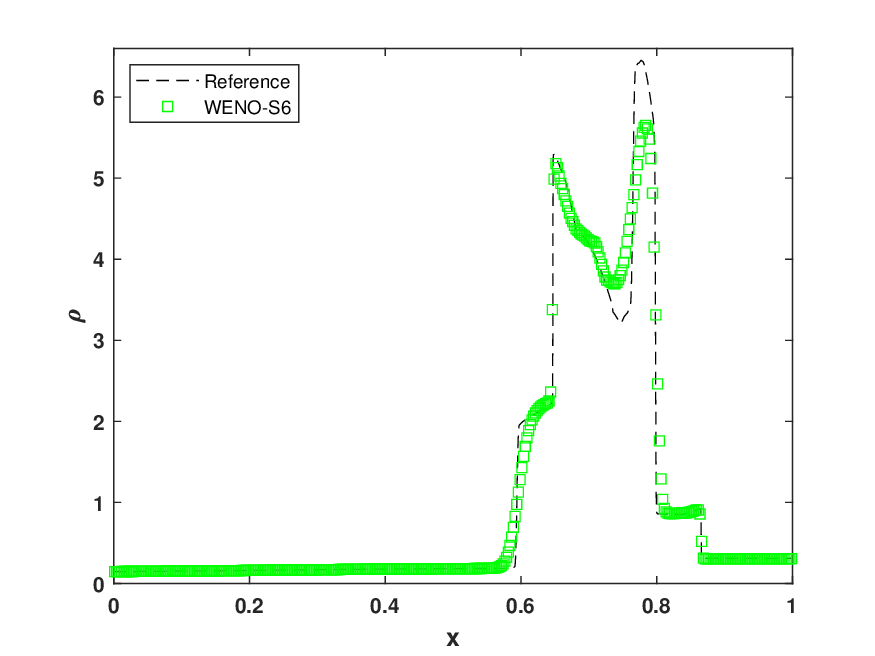}
\includegraphics[width=0.245\textwidth, trim=30 5 30 10, clip=true]{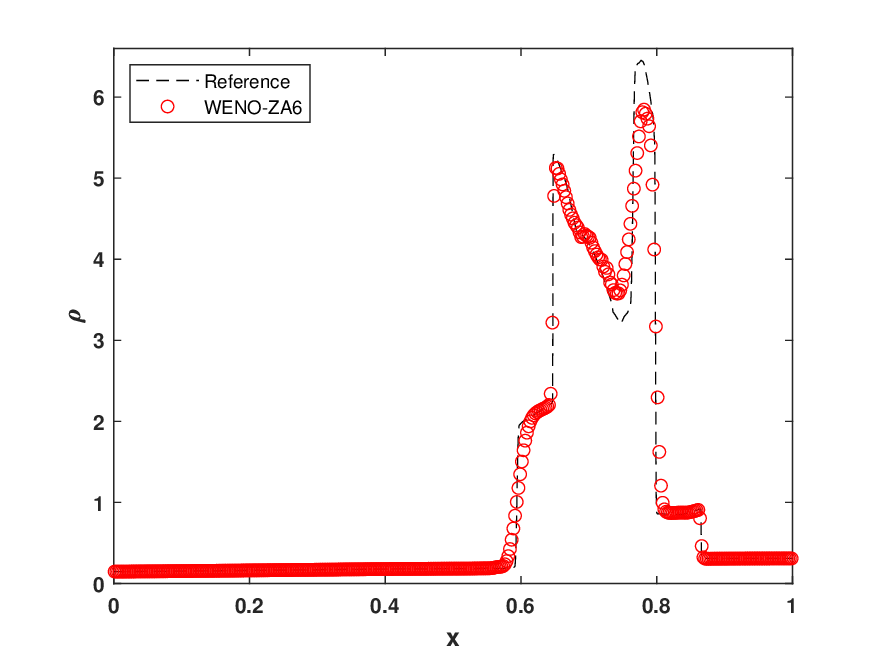}
\includegraphics[width=0.245\textwidth, trim=30 5 30 10, clip=true]{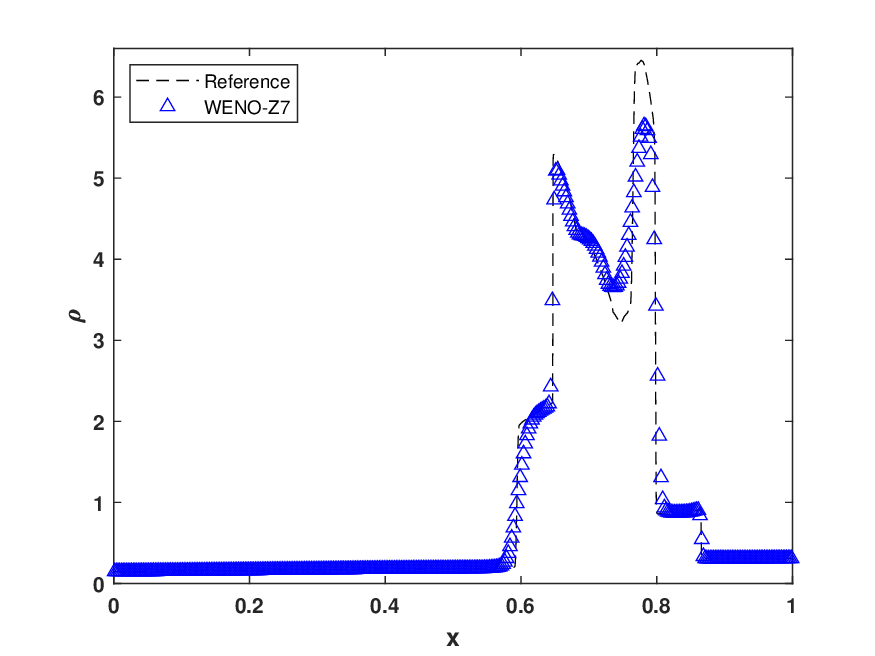}
}
\end{center}
\vskip-10pt
\caption{Example~\ref{ex:blast}, (Left) numerical density of the blastwaves interaction problem and (Right) close-up view of the density by the WENO-CU6 (purple), WENO-S6 (green), WENO-ZA6 (red), and WENO-Z7 (blue) schemes with $N = 400$ at $T = 0.038$.
}
\label{fig:blast}
\end{figure}

\begin{example}\label{ex:blast}\rm 
We consider the blastwaves interaction problem \cite{Woodward} with the initial condition
$$ 
\left(\rho,~u,~P\right) = 
\left\{
\begin{array}{llllcccc}
(& 1, & 0, & 1000 &), & 0.0 \leq x < 0.1,\\
(& 1, & 0, & 0.01 &), & 0.1 \leq x < 0.9,\\
(& 1, & 0, & 100  &), & 0.9 \leq x \leq 1.0.
\end{array}
\right.
$$
The computational domain $[0,1]$ is divided into $N=400$ uniform cells.
The numerical solutions for the density \(\rho\) at the final time \(T = 0.038\) are displayed in Fig.~\ref{fig:blast}. 
The reference solution is computed using WENO-Z7 on a refined grid with \(N = 4000\) cells. 
It can be observed that all numerical solutions are in good agreement with the reference profile. 
Furthermore, the enlarged view indicates that WENO-ZA6 achieves slightly higher resolution than the other three WENO schemes.
\end{example}

\subsection{Two-dimensional Euler equations} \label{sec:Euler2d}

In this section, we solve the two-dimensional Euler equations and present numerical results for four benchmark problems: the Riemann initial-value problem (RIVP), the double Mach reflection problem, the single-material triple-point problem, and the shock diffraction problem. These test cases involve complex flow features, including strong gradients, discontinuities, and small-amplitude smooth structures, thereby providing a stringent assessment of the resolution and robustness of the proposed WENO-ZA6 scheme.

\begin{figure}[htb!]
\begin{center}
\mbox{
\makebox[0.4\textwidth][c]{\hskip10pt WENO-CU6}
\makebox[0.4\textwidth][c]{\hskip10pt WENO-S6}
}
\mbox{
\includegraphics[width=0.4\textwidth, trim=50 10 70 10, clip=true]{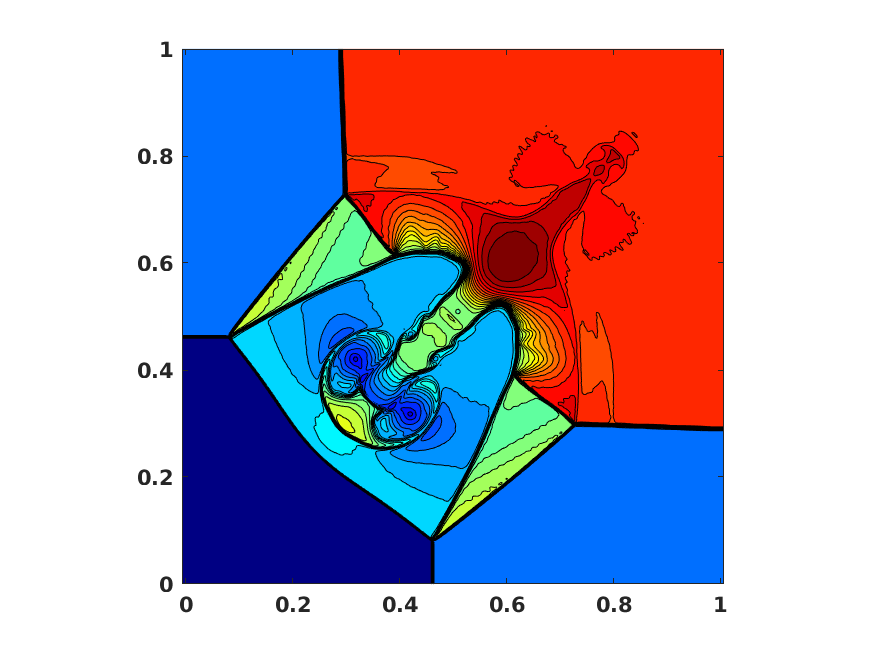}
\includegraphics[width=0.4\textwidth, trim=50 10 70 10, clip=true]{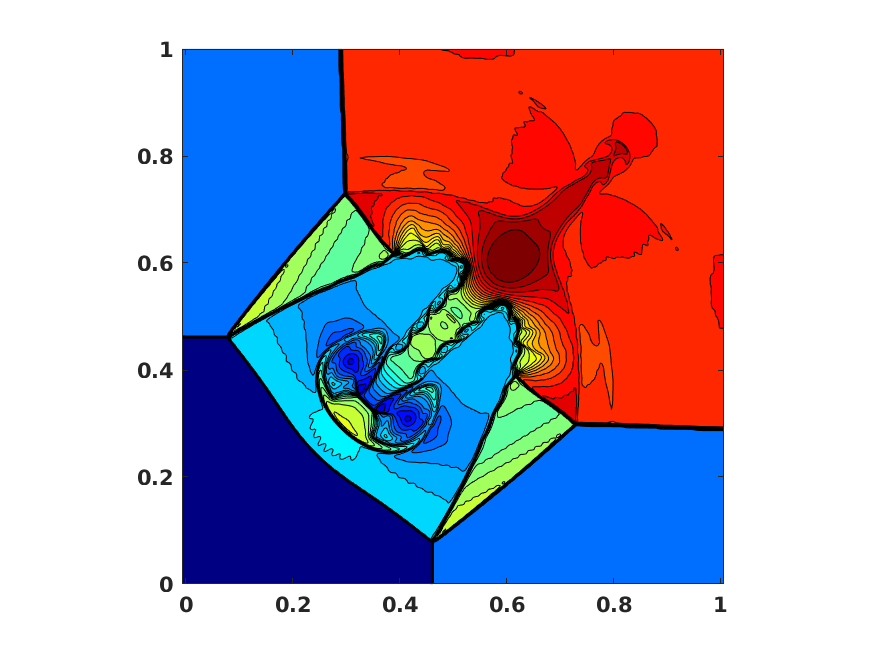}
}
\mbox{
\makebox[0.4\textwidth][c]{\hskip10pt WENO-ZA6}
\makebox[0.4\textwidth][c]{\hskip10pt WENO-Z7}
}
\mbox{
\includegraphics[width=0.4\textwidth, trim=50 10 70 10, clip=true]{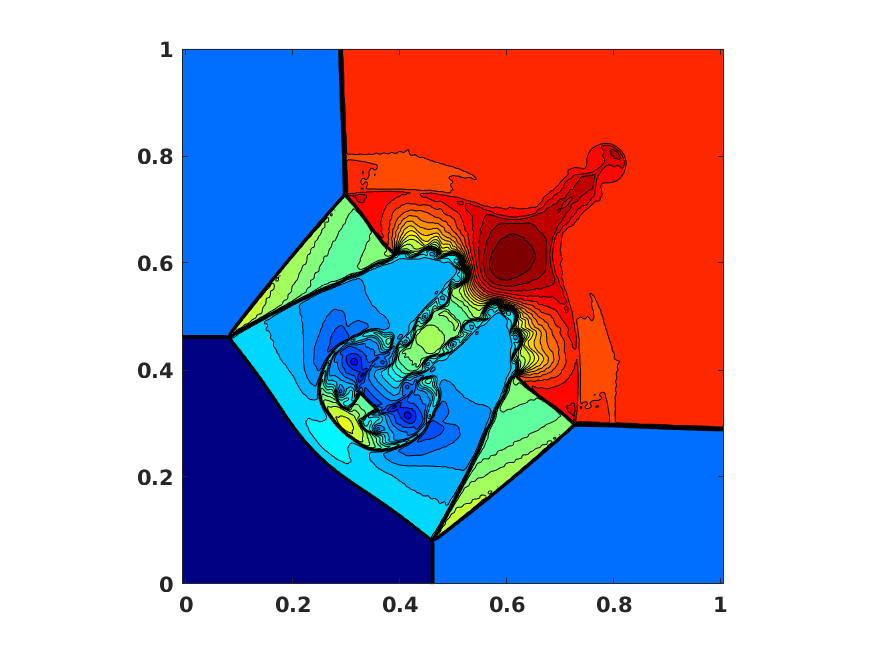}
\includegraphics[width=0.4\textwidth, trim=50 10 70 10, clip=true]{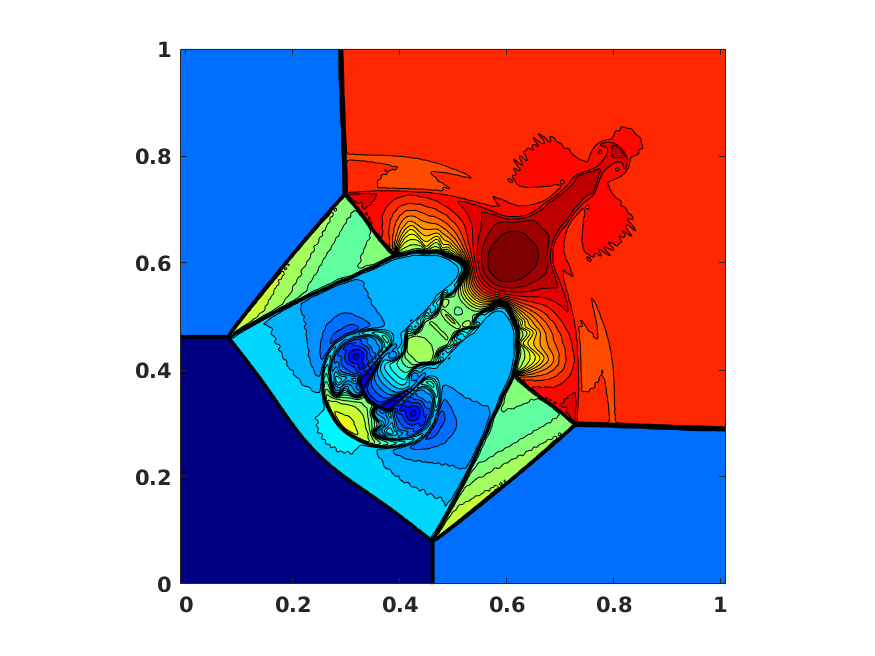}
}
\end{center}
\vskip-10pt
\caption{Example \ref{ex:Riemann_config3}, the density of configuration 3 of the Riemann initial value problem computed by the WENO schemes with $N_x = N_y = 400$ at $T = 0.8$.}
\label{fig:Riemann_config3}
\end{figure}

\begin{example}\label{ex:Riemann_config3}\rm 
Configuration 3 of the two-dimensional Riemann problem \cite{LaxLiu} is tested here. 
The initial condition contains four quadrants of constant states separated by lines $x = 0.8$ and $y = 0.8$,
$$ 
\left(\rho,~u,~v,~P\right) = 
\left\{
\begin{array}{lllllccc}
(& 1.5,    & 0,     & 0,     & 1.5   &), & x>0.8, \: y>0.8, \\
(& 0.5323, & 1.206, & 0,     & 0.3   &), & x<0.8, \: y>0.8, \\
(& 0.138,  & 1.206, & 1.206, & 0.029 &), & x<0.8, \: y<0.8, \\
(& 0.5323, & 0,     & 1.206, & 0.3   &), & x>0.8, \: y<0.8.
\end{array}
\right.
$$
We divide the square domain $[0, 1] \times [0, 1]$ into $N_x \times N_y = 400 \times 400$ uniform cells.
The simulations of the density computed by the WENO schemes at the final time, $T = 0.8$, are displayed in Fig. \ref{fig:Riemann_config3}. 
The results show that the flow's large-scale structures, including the triple point, the incident shock, the reflected shock, a Mach stem, and a slip plane, are in good agreement with those reported in \cite{DonLiWangWang}. 
\end{example}

\begin{figure}[htb!]
\begin{center}
\mbox{
\makebox[0.45\textwidth][c]{WENO-CU6}
\makebox[0.45\textwidth][c]{WENO-S6 }
}
\mbox{
\includegraphics[width=0.45\textwidth, trim=30 10 30 5, clip=true]{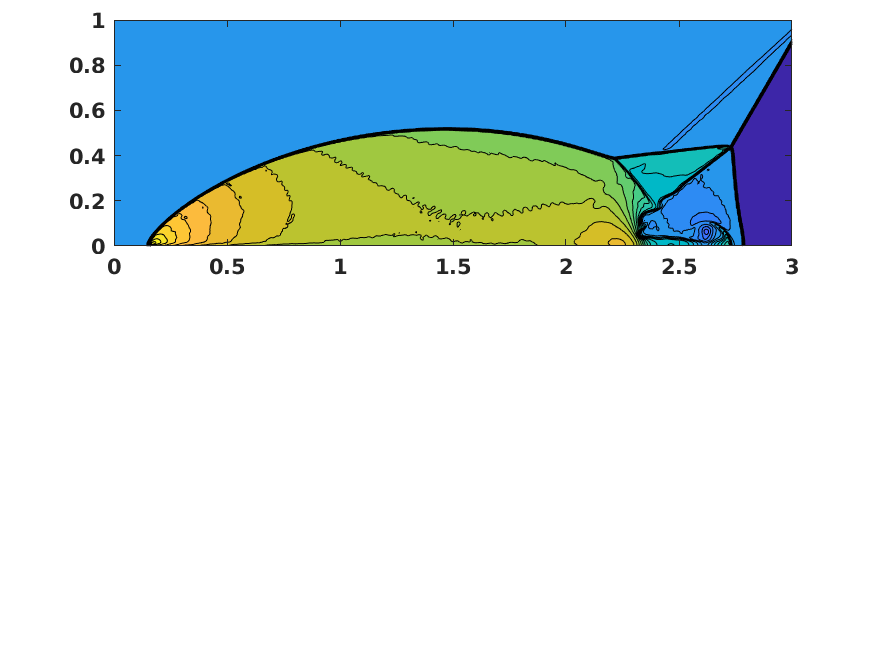}
\includegraphics[width=0.45\textwidth, trim=30 10 30 5, clip=true]{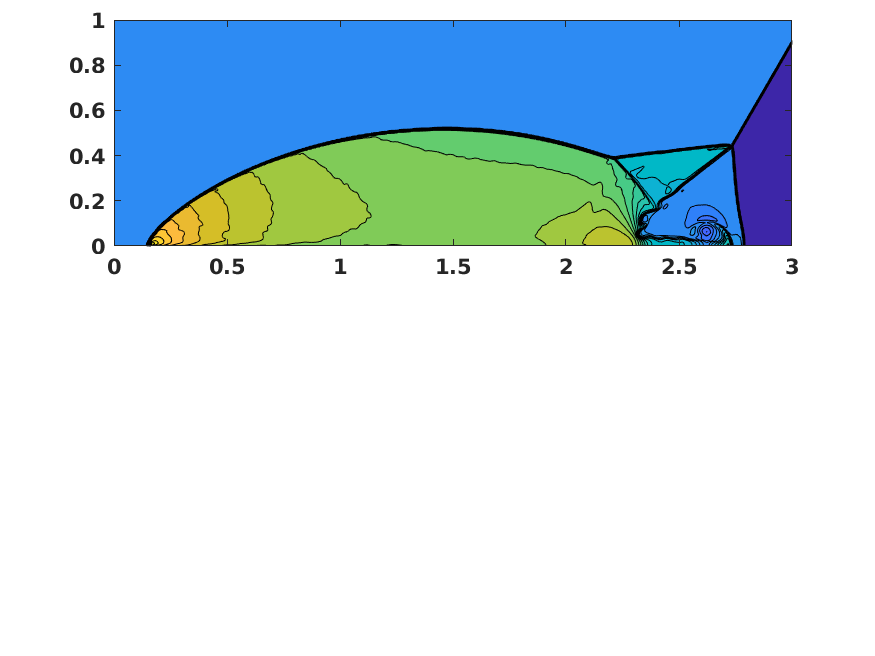}
}
\vskip-108pt
\mbox{
\makebox[0.45\textwidth][c]{WENO-ZA6}
\makebox[0.45\textwidth][c]{WENO-Z7}
}
\mbox{
\includegraphics[width=0.45\textwidth, trim=30 10 30 5, clip=true]{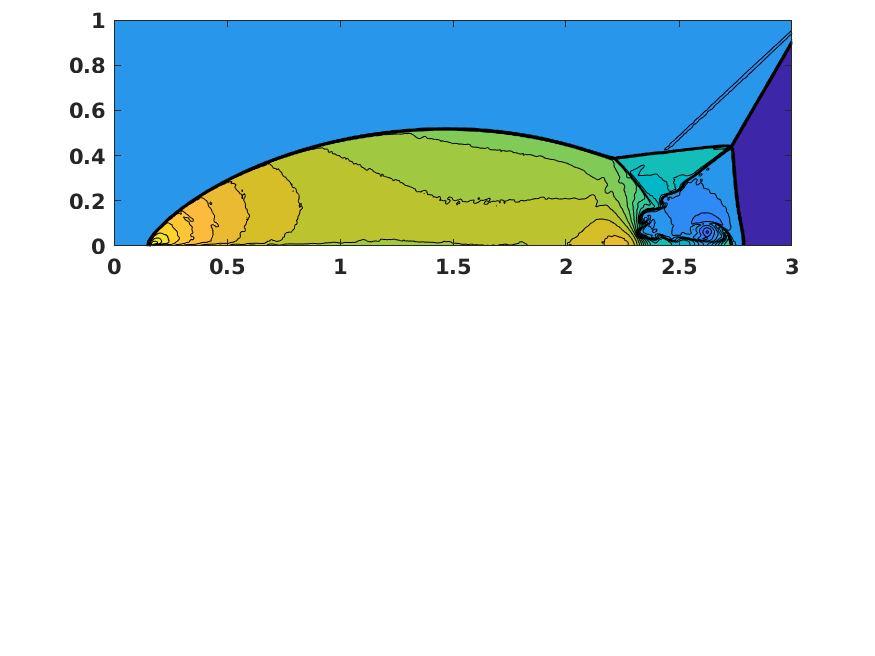}
\includegraphics[width=0.45\textwidth, trim=30 10 30 5, clip=true]{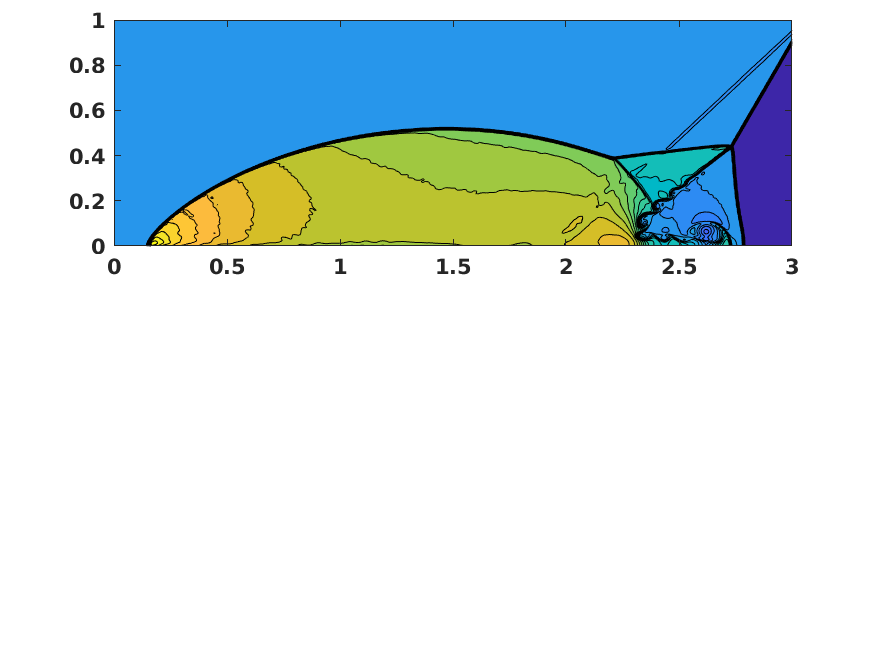}
}
\end{center}
\vskip-120pt
\caption{Example~\ref{ex:DMR}, the density of double Mach reflection problem computed by the WENO schemes with $N_x \times N_y = 800 \times 200$ at $T = 0.2$.}
\label{fig:DMR}
\end{figure}

\begin{figure}[htb!]
\begin{center}
\mbox{
\makebox[0.45\textwidth][c]{WENO-CU6}
\makebox[0.45\textwidth][c]{WENO-S6}
}
\mbox{
\includegraphics[width=0.45\textwidth, trim=30 10 30 0, clip=true]{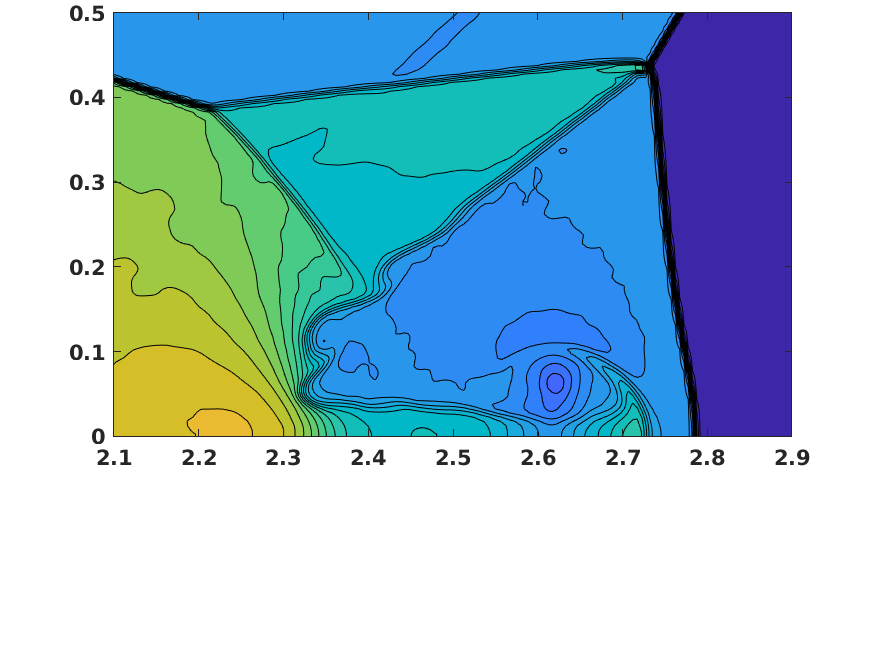}
\includegraphics[width=0.45\textwidth, trim=30 10 30 0, clip=true]{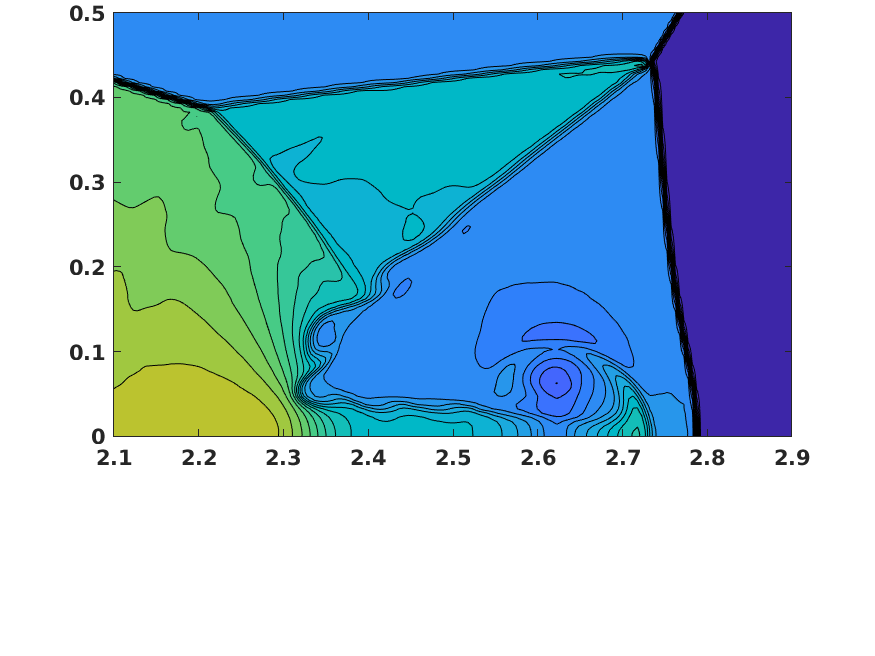}
}
\vskip-50pt
\mbox{
\makebox[0.45\textwidth][c]{WENO-ZA6}
\makebox[0.45\textwidth][c]{WENO-Z7}
}
\mbox{
\includegraphics[width=0.45\textwidth, trim=30 10 30 0, clip=true]{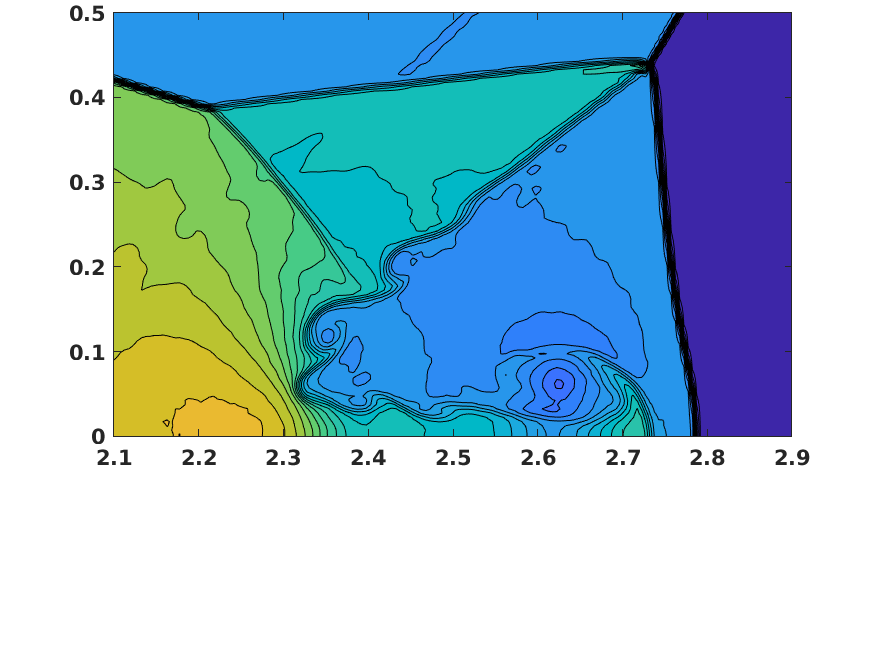}
\includegraphics[width=0.45\textwidth, trim=30 10 30 0, clip=true]{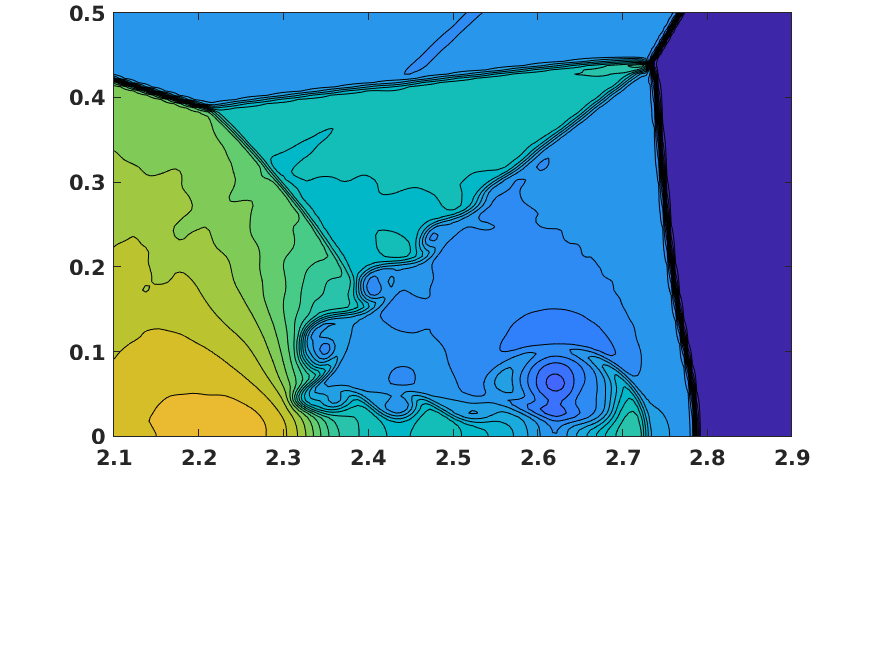}
}
\end{center}
\vskip-65pt
\caption{Example \ref{ex:DMR}, close-up view of the density of double Mach reflection problem computed by the WENO schemes with $N_x \times N_y = 800 \times 200$ at $T = 0.2$.}
\label{fig:DMR_zoomin}
\end{figure}

\begin{example}\label{ex:DMR}\rm 
We consider the two-dimensional double Mach reflection of a shock from an oblique surface. 
It describes the reflection of a planar Mach shock in the air striking a wedge, as proposed and studied in detail by Woodward and Colella \cite{Woodward}; the wall jet curls into a vortex as the jet reaches the Mach stem. 
Furthermore, the shock that joins the contact surface and transverse wave (making this a double Mach reflection) sharpens, and its triple point, the kink along the transverse wave, becomes much more distinguished. 
Finally, Kelvin-Helmholtz instabilities develop along the wall jet's contact surface.

It is computed on $[0,4] \times [0,1]$ and the initial condition is
$$
\left(\rho,~u,~v,~P\right) = 
\left\{
\begin{array}{cccccccr}
(& 8.0, & 8.25 \,\cos(\pi/6), & -8.25 \, \sin(\pi/6), & 116.5 &), & x < \frac{1}{6} + \frac{y}{\sqrt{3}} , \\
(& 1.4, & 0,                  & 0,                    &  1    &), & x \geq \frac{1}{6} + \frac{y}{\sqrt{3}}.
\end{array}
\right.
$$
We use the absolute value in the computation of the sound speed for WENO-S6 to prevent rapid blow-up of the numerical solution.

Figure \ref{fig:DMR} shows the density profiles computed by the WENO schemes at the final time $T = 0.2$ on the uniform mesh with $N_x \times N_y = 800 \times 200$ cells.
We further zoom in on the solution for the region $[2.1, 2.9] \times [0, 0.5]$ in Fig. \ref{fig:DMR_zoomin}.
The results show that WENO-ZA6 and WENO-Z7 effectively capture shocks and high gradients. 
Additionally, these schemes are more effective than WENO-CU6 and WENO-S6 at resolving small, weak, high-frequency structures along the slip line and the vortical rollup at the jet tip. 
\end{example}

\begin{figure}[htb!]
\begin{center}
\mbox{
\makebox[0.45\textwidth][c]{WENO-CU6}
\makebox[0.45\textwidth][c]{WENO-S6}
}
\mbox{
\includegraphics[width=0.45\textwidth, trim=30 10 30 5, clip=true]{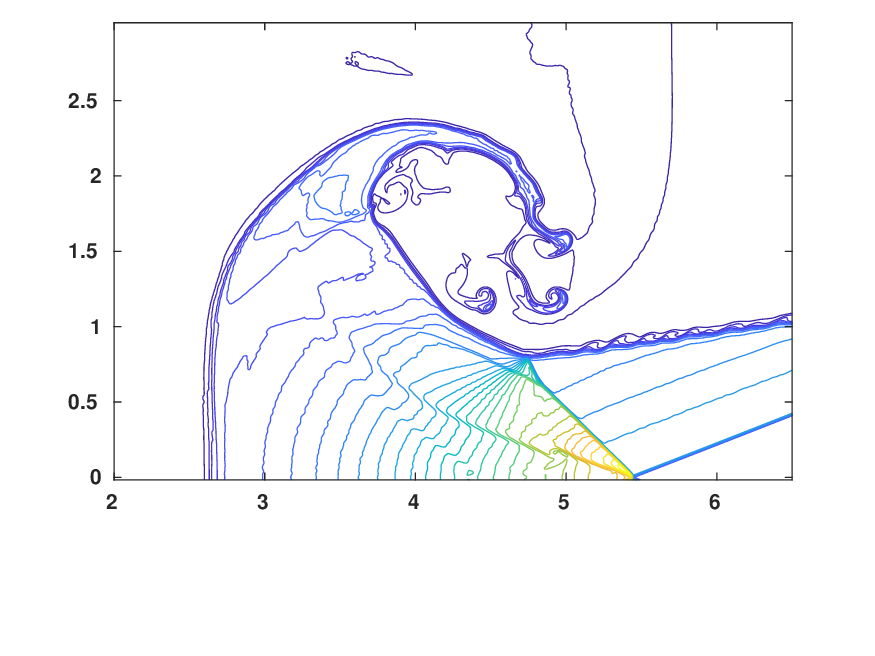}
\includegraphics[width=0.45\textwidth, trim=30 10 30 5, clip=true]{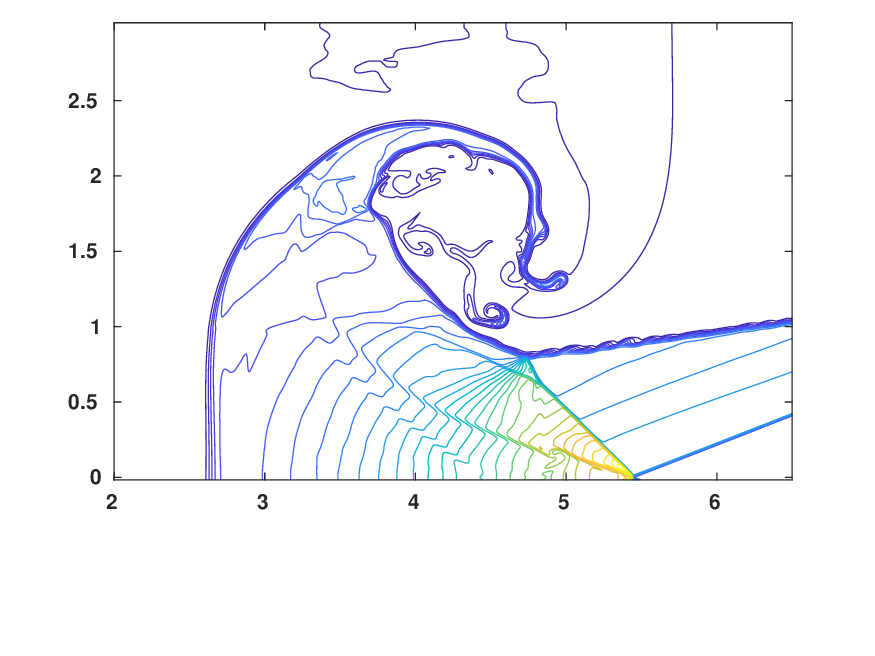}
}
\vskip-39pt
\mbox{
\makebox[0.45\textwidth][c]{WENO-ZA6}
\makebox[0.45\textwidth][c]{WENO-Z7}
}
\mbox{
\includegraphics[width=0.45\textwidth, trim=30 10 30 5, clip=true]{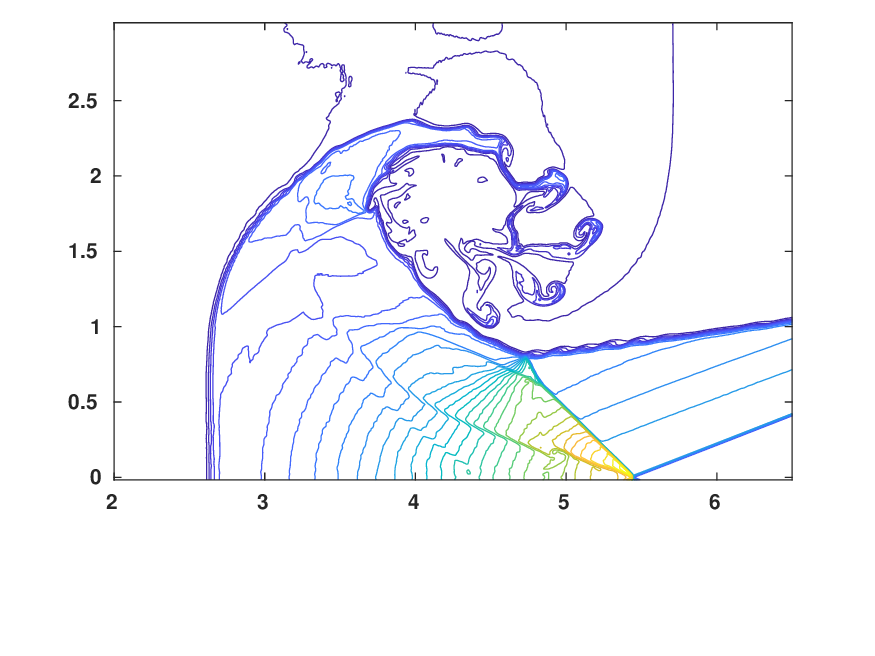}
\includegraphics[width=0.45\textwidth, trim=30 10 30 5, clip=true]{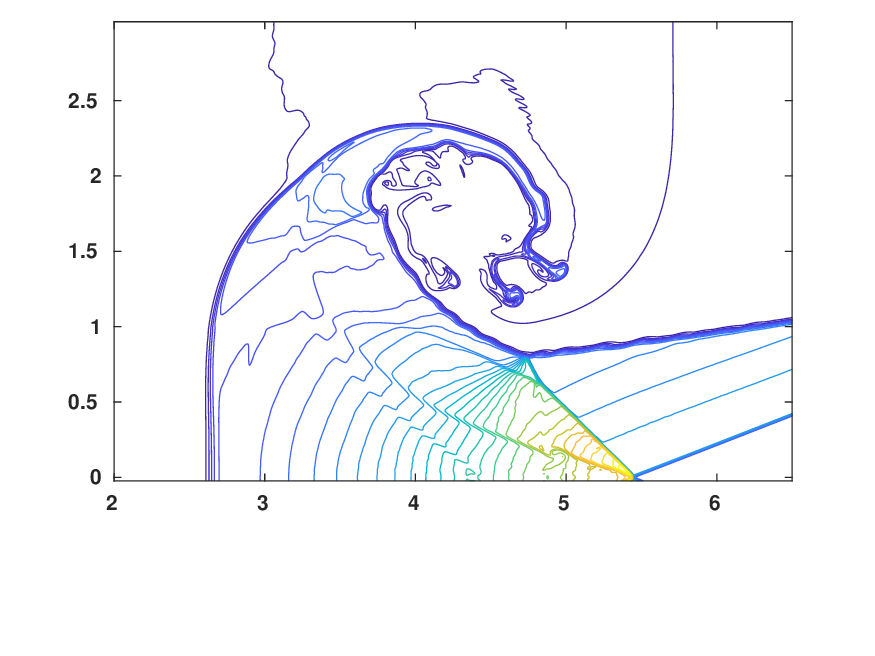}
}
\end{center}
\vskip-50pt
\caption{Example \ref{ex:SMTP}, the density of single-material triple point problem computed by the WENO schemes with $N_x \times N_y = 1050 \times 450$ at $T = 5.0$.}
\label{fig:SMTP}
\end{figure}

\begin{example}\label{ex:SMTP}\rm 
In this example, we investigate the modified triple-point problem (SMTP) \cite{Zeng} using a single material. 
The initial condition is given by
$$ 
\left(\rho,~u,~v,~P\right) = 
\left\{
\begin{array}{llllccl}
(& 1.000, & 0, & 0, & 1.0 &), & x < 1.0 , \\
(& 1.000, & 0, & 0, & 0.1 &), & x >1.0, ~y<1.5, \\
(& 0.125, & 0, & 0, & 0.1 &), & x >1.0, ~y>1.5.
\end{array}
\right.
$$
The computational domain is $[0,7] \times [0,3]$ with the uniform $N_x \times N_y = 1050 \times 450$ cells, and the final time is set to $T=5$. 
Fig. \ref{fig:SMTP} shows the density computed by the WENO schemes at the final time, indicating that the large-scale structures produced by all WENO schemes are similar. 
Besides, the small-scale structures along the slip line are much better resolved by WENO-ZA6 due to its higher resolution and lower dissipation.
\end{example}

\begin{figure}[htb!]
\begin{center}
\mbox{
\makebox[0.4\textwidth][c]{WENO-CU6}\quad
\makebox[0.4\textwidth][c]{WENO-S6}
}
\mbox{
\includegraphics[width=0.4\textwidth, trim=50 10 50 20, clip=true]{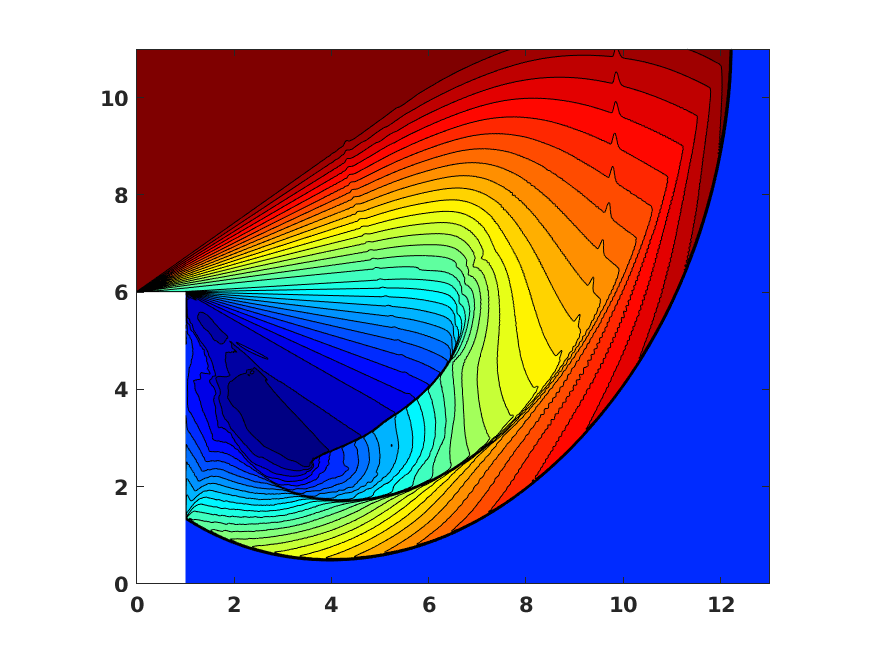}\quad
\includegraphics[width=0.4\textwidth, trim=50 10 50 20, clip=true]{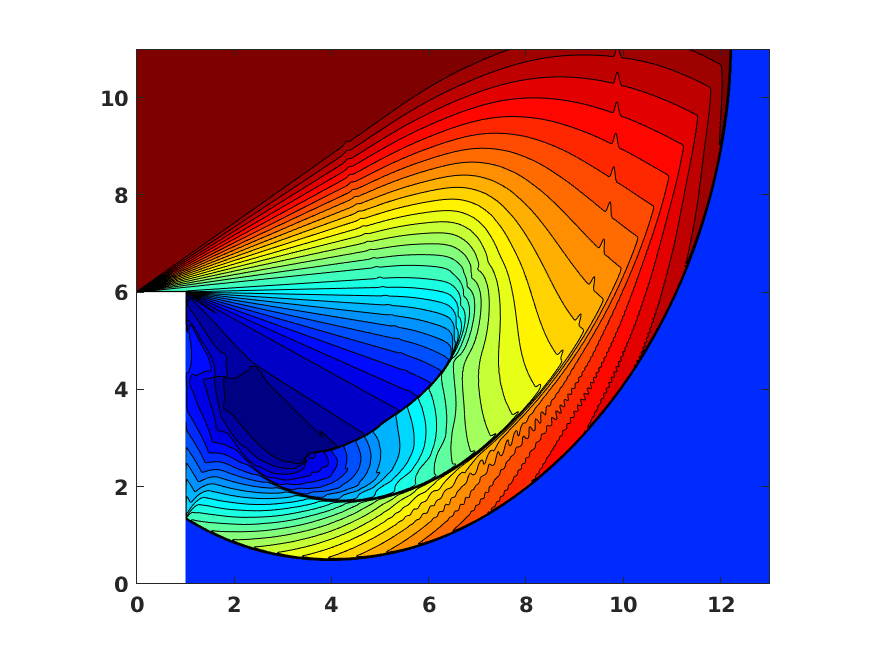}
}
\mbox{
\makebox[0.4\textwidth][c]{WENO-ZA6}\quad
\makebox[0.4\textwidth][c]{WENO-Z7}
}
\mbox{
\includegraphics[width=0.4\textwidth, trim=50 10 50 20, clip=true]{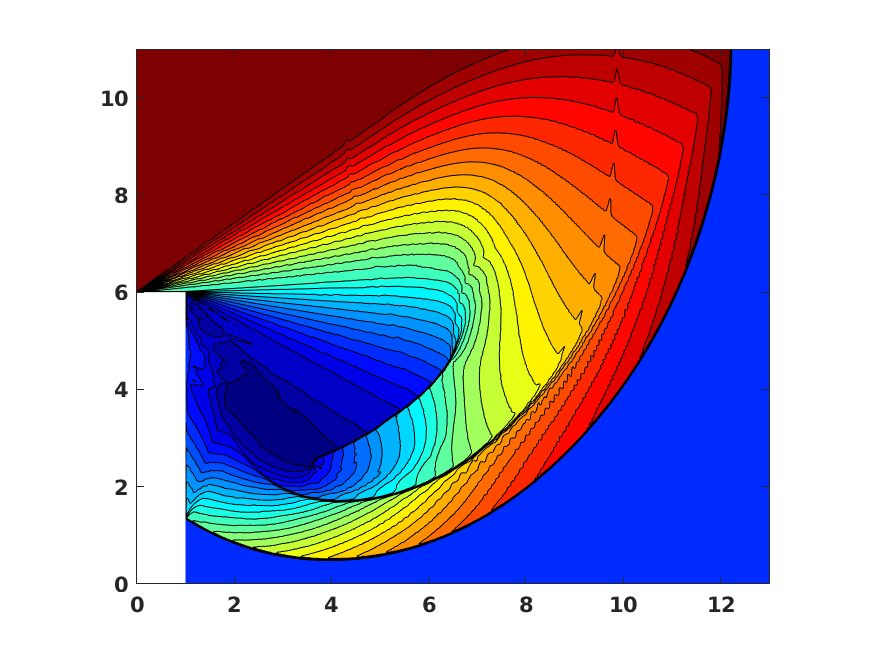}\quad
\includegraphics[width=0.4\textwidth, trim=50 10 50 20, clip=true]{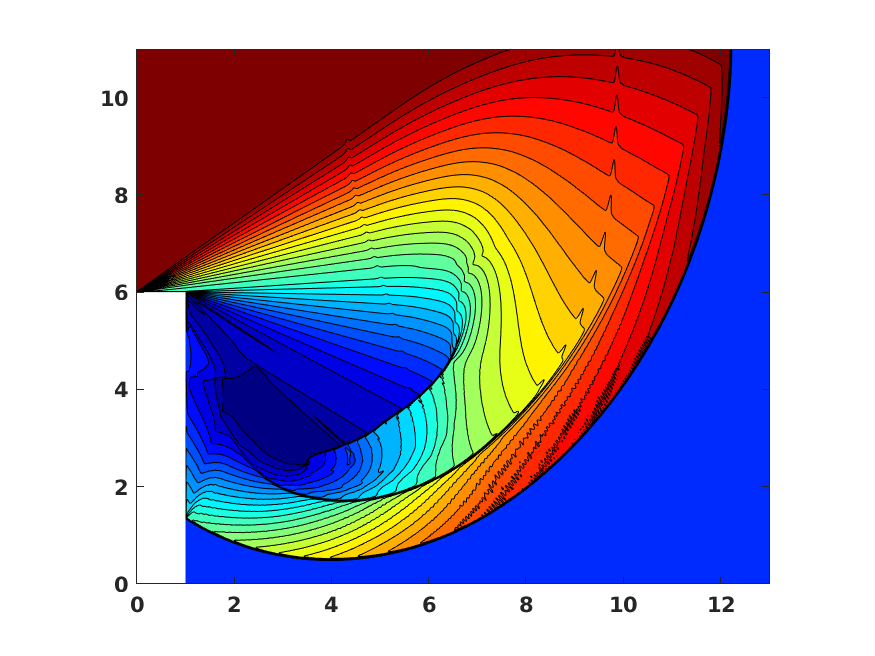}
}
\end{center}
\vskip-10pt
\caption{Example \ref{ex:SD}, the density of shock diffraction problem by the WENO schemes with the positivity-preserving limiter for $N_x \times N_y = 1040 \times 880$ at $T = 2.3$.}
\label{fig:SD}
\end{figure}

\begin{example}\label{ex:SD}\rm 
We end this subsection with the shock diffraction problem \cite{CockburnShu, ZhangShu10} to examine the robustness of the WENO schemes with the positivity-preserving limiter. 
The computational domain is the union of $[0,1] \times [6,11]$ and $[1,13] \times [0,11]$. 
The initial condition is a pure right-moving shock with Mach number $5.09$, located at $x = 0.5$ and $6 \leq y \leq 11$, moving into undisturbed air ahead of it, with density $\rho=1.4$ and pressure $P=1$. 
We take $\gamma = 1.4$.
The boundary conditions are inflow at $x=0, 6\leq y \leq11$, outflow at $x=13, 0\leq y \leq11$, $y=11, 0\leq x \leq13$ and $y=0, 1\leq x \leq13$, and reflective at the walls $0\leq x \leq1, y=6$ and $x=1, 0\leq y \leq6$. 
The density contours at the final time $T = 2.3$ computed with $N_x \times N_y = 1040 \times 880$ are presented in Fig. \ref{fig:SD}.

The shock diffraction problem is numerically challenging, especially for high-order schemes, because the pressure and density may approach zero when the shock wave diffracts around the obstacle, leading to a near-vacuum state. 
To handle this extreme situation, the positivity-preserving limiter~\cite{HuAdamsShu,LiDonWangGao,LiDonWangWang} is applied to maintain positivity in troubled cells, which typically appear below and to the right of the corner at early times. 
The numerical results show that the large-scale flow structures produced by the different WENO schemes are essentially identical. 
Owing to its less dissipation, WENO-ZA6 significantly resolves small-scale structures along the slip line and improves vortex-capturing capabilities.
\end{example}  

\subsection{CPU timings of WENO schemes}
Before concluding, we compare the computational efficiency of the WENO-CU6, WENO-S6, WENO-ZA6, and WENO-Z7, using the two-dimensional Euler equations considered in the previous subsection. 
The Fortran~90 code is compiled with the INTEL \texttt{ifort} compiler using the \texttt{-O2} optimization flag. All simulations are carried out on an eight-core Intel i7-9700 CPU running at 3.00\, GHz. Since the relative efficiency is nearly identical across all two-dimensional test cases, we report only the results for Example~\ref{ex:Riemann_config3} as a representative example.

Table~\ref{tab:CPU_comparison} presents the CPU time together with the corresponding speedup relative to WENO-CU6 as the reference, defined as
$ \text{speedup} = \frac{\text{CPU}_{\text{scheme}}-\text{CPU}_{\text{CU6}}}{\text{CPU}_{\text{CU6}}} \times 100\%$.
Among the four schemes, WENO-ZA6 requires the least CPU time per Runge--Kutta step, primarily due to its simpler construction and formulation of the nonlinear weights. 
Overall, the proposed WENO-ZA6 is comparable in cost to WENO-CU6, approximately 21\% faster than WENO-S6, and about 15\% faster than WENO-Z7.

\begin{table}[htb!]
\begin{center}
\caption{CPU times and speedup for Example \ref{ex:Riemann_config3}.}
\label{tab:CPU_comparison}
\begin{tabular}{l c c c c}
\hline
& WENO-CU6 & WENO-S6 & WENO-ZA6 & WENO-Z7 \\
\hline
CPU time          & 2722 & 3313  & 2608 & 3174 \\
\# of steps       & 3495 & 3522  & 3497 & 3532 \\
Speedup (\%)      & --   & -21\% & 4\%  & -15\%  \\
\hline
\end{tabular}
\end{center}
\end{table}

\section{Conclusion} 
\label{sec:conclusion}

In this paper, we developed a systematic WENO reconstruction framework for even-order central-upwind WENO finite-difference schemes using a simple downwind smoothness indicator defined as the arithmetic mean of local smoothness indicators. This approach incorporates regularity information across the stencil without added complexity or the need for tuning parameters. Analysis of the sixth-order scheme (WENO-ZA6) shows that the nonlinear weights ensure optimal convergence at critical points. Dispersion analysis confirms WENO-ZA6 mirrors the underlying linear scheme's characteristics while maintaining essential numerical dissipation for robust shock capturing.

Comprehensive numerical experiments support these theoretical findings. For scalar conservation laws with both convex and nonconvex fluxes, WENO-ZA6 accurately captures isolated shocks and compound wave structures in a non-oscillatory manner. In the context of the one-dimensional Euler equations, the scheme solves several classical Riemann problems, shock-entropy interactions, extended shock-density waves, and the two-blast-wave interaction problem well. Two-dimensional tests, such as the Riemann initial-value problem, double Mach reflection, triple-point problem, and shock diffraction, further demonstrate the scheme's ability to capture complex flow features with high resolution and essentially non-oscillatory shock capturing. Notably, WENO-ZA6 achieves the same computational efficiency as WENO-CU6 and provides substantial gains of approximately 15–21\% over WENO-S6 and WENO-Z7 due to its simpler weight formulation.

Future work includes applying this framework to finite-volume schemes and to more complex systems such as steady-state compressible flows, magnetohydrodynamics, and multi-component flows.

\appendix
\renewcommand{\appendixname}{Appendix~\Alph{section}}
\section*{Appendix}
\section{Fourth-, eighth-, and tenth-order central-upwind WENO schemes}
\label{sec:appendix}

In this section, we provide the explicit coefficient vectors and matrices corresponding to the fourth-, eighth-, and tenth-order central-upwind WENO schemes for reference. Furthermore, a symbolic MATLAB program is included to allow readers to generate and validate the coefficients for central-upwind WENO schemes of arbitrary even order.

\subsection{Fourth-order central-upwind WENO scheme}

For the fourth-order central-upwind WENO scheme with the positive advection speed (or the $+$ split flux), we consider the $4$-point stencil $S^4$, consisting of two upwind substencils $\{ S_0, S_1 \}$ and one downwind substencil $S_3$. 
For the stencil $S^4$, the coefficient vector in \eqref{eq:numerical_flux_FD_plus} is
$$ \bfc_4 = \tfrac{1}{12}(-1,7,7,-1). $$
The coefficient vectors $\bfc_k$ in \eqref{eq:numerical_subflux} associated with the $2$-point substencils $S_k$ are
\begin{equation*}\label{eq:subflux_S4_coeffs}
\left\{
\begin{aligned}
\bfc_0 &= \tfrac{1}{2}(-1,3), \\
\bfc_1 &= \tfrac{1}{2}(1,1), 
\end{aligned}
\right.
\qquad \text{and} \quad
\bfc_2 = \mathrm{rev}(\bfc_0).
\end{equation*}
The ideal weight vector in \eqref{eq:numerical_flux_FD_linear_comb} is
$$ \bfd = \tfrac{1}{6}(1,4,1). $$
 Using \eqref{eq:SI}, the smoothness indicators $\beta_k$ can be written in the unified form
$$
\beta_k = \langle \bff_k, \mathbf{B}_k \bff_k \rangle, \quad k \in \{0,1,d\},
$$
where the corresponding smoothness matrices are 
\begin{eqnarray*}
\mathbf{B}_0 &=&
\left[
\begin{array}{rr}
  1 & -1 \\
 -1 &  1
\end{array}
\right],
\qquad
\mathbf{B}_d =
\left[
\begin{array}{rr}
  1 & -1 \\
 -1 &  1
\end{array}
\right],
\end{eqnarray*}
and $\mathbf{B}_{1} = \mathrm{rev}(\mathbf{B}_0)$.
The smoothness indicator for the downwind substencil $S_2$ is then defined by 
$$ \beta_2= \frac{1}{3} \left( \beta_0 + \beta_1 + \beta_d \right). $$
Finally, the global smoothness indicator $\tau$ \eqref{eq:GSI} employs the anti-symmetric vector 
$$ \bfc_\tau = (-1,3,-3,1). $$

\subsection{Eighth-order central-upwind WENO scheme}

For the eighth-order central-upwind WENO scheme with the positive advection speed (or the $+$ split flux), we consider the $8$-point stencil $S^8$, consisting of four upwind substencils $\{ S_0, S_1, S_2, S_3 \}$ and one downwind substencil $S_4$, as shown in Fig. \ref{fig:stencil8}.
For the stencil $S^8$, the coefficient vector in \eqref{eq:numerical_flux_FD_plus} is
$$ \bfc_8 = \tfrac{1}{840}(-3,29,-139,533,533,-139,29,-3). $$
The coefficient vectors $\bfc_k$ in \eqref{eq:numerical_subflux} associated with the $4$-point substencils $S_k$ are
\begin{equation*}\label{eq:subflux_S8_coeffs}
\left\{
\begin{aligned}
\bfc_0 &= \tfrac{1}{12}(-3,13,-23,25), \\
\bfc_1 &= \tfrac{1}{12}(1,-5,13,3), \\
\bfc_2 &= \tfrac{1}{12}(-1,7,7,-1),
\end{aligned}
\right.
\qquad \text{and} \quad
\bfc_{4-k} = \mathrm{rev}(\bfc_k), 
\ k=0,1.
\end{equation*}
The ideal weight vector in \eqref{eq:numerical_flux_FD_linear_comb} is
$$ \bfd = \tfrac{1}{70}(1,16,36,16,1). $$
Using \eqref{eq:SI}, the smoothness indicators $\beta_k$ can be written in the unified form
$$
\beta_k = \langle \bff_k, \mathbf{B}_k \bff_k \rangle, \quad k \in \{0,1,2,3,d\},
$$
where the corresponding smoothness matrices are 
\begin{eqnarray*}
\mathbf{B}_0 &=& \scriptsize \tfrac{1}{240}
\left[
\begin{array}{rrrr}
   547 & -1941 &  2321 &  -927 \\
 -1941 &  7043 & -8623 &  3521 \\
  2321 & -8623 & 11003 & -4701 \\
  -927 &  3521 & -4701 &  2107
\end{array}
\right],
\qquad 
\mathbf{B}_1 = \scriptsize \tfrac{1}{240}
\left[
\begin{array}{rrrr}
   267 &  -821 &   801 &  -247 \\
  -821 &  2843 & -2983 &   961 \\
   801 & -2983 &  3443 & -1261 \\
  -247 &   961 & -1261 &   547
\end{array}
\right],
\\
\mathbf{B}_d &=& \scriptsize \tfrac{1}{240}
\left[
\begin{array}{rrrr}
   7107 & -16901 &  13521 & -3727 \\
 -16901 &  40643 & -32863 &  9121 \\
  13521 & -32863 &  26843 & -7501 \\
  -3727 &   9121 &  -7501 &  2107
\end{array}
\right],
\end{eqnarray*}
and $\mathbf{B}_{3-k} = \mathrm{rev}(\mathbf{B}_k), k = 0,1$.
The smoothness indicator for the downwind substencil $S_4$ is then defined by 
$$ \beta_4= \frac{1}{5} \left( \sum_{k=0}^3 \beta_k + \beta_d \right). $$
Finally, the global smoothness indicator $\tau$ \eqref{eq:GSI} employs the anti-symmetric vector 
$$ \bfc_\tau = (-1,7,-21,35,-35,21,-7,1). $$

\subsection{Tenth-order central-upwind WENO scheme}

For the tenth-order central-upwind WENO scheme with the positive advection speed (or the $+$ split flux), we use the $10$-point stencil $S^{10}$, consisting of five upwind substencils $\{ S_0, S_1, S_2, S_3, S_4 \}$ and one downwind substencil $S_5$.
For the stencil $S^{10}$, the coefficient vector in \eqref{eq:numerical_flux_FD_plus} is
$$ \bfc_{10}= \tfrac{1}{2520}(2,-23,127,-473,1627,1627,-473,127,-23,2). $$
The coefficient vectors $\bfc_k$ in \eqref{eq:numerical_subflux} associated with the $5$-point substencils $S_k$ are
\begin{equation*}\label{eq:subflux_S10_coeffs}
\left\{
\begin{aligned}
\bfc_0 &= \tfrac{1}{60}(12,-63,137,-163,137), \\
\bfc_1 &= \tfrac{1}{60}(-3,17,-43,77,12), \\
\bfc_2 &= \tfrac{1}{60}(2,-13,47,27,-3), 
\end{aligned}
\right.
\qquad \text{and} \quad \bfc_{5-k} = \mathrm{rev}(\bfc_k), \ k=0,1,2.
\end{equation*}
The ideal weight vector in \eqref{eq:numerical_flux_FD_linear_comb} is
$$ \bfd = \tfrac{1}{252}(1,25,100,100,25,1). $$
Using \eqref{eq:SI}, the smoothness indicators $\beta_k$ can be written in the unified form
$$
\beta_k = \langle \bff_k, \mathbf{B}_k \bff_k \rangle, \quad k \in \{0,1,2,3,4,d\},
$$
where the corresponding smoothness matrices are 
\begin{eqnarray*} 
\mathbf{B}_0 &=& \scriptsize \tfrac{1}{10080}
\left[ 
\begin{array}{rrrrr}
   45316 &  -208501 &   364863 &  -288007 &   86329 \\
 -208501 &   965926 & -1704396 &  1358458 & -411487 \\
  364863 & -1704396 &  3042786 & -2462076 &  758823 \\
 -288007 &  1358458 & -2462076 &  2041126 & -649501 \\
   86329 &  -411487 &   758823 &  -649501 &  215836
\end{array}
\right],
\\
\mathbf{B}_1 &=& \scriptsize \tfrac{1}{10080}
\left[
\begin{array}{rrrrr}
   13816 &  -60871 &   99213 &  -70237 &   18079 \\
  -60871 &  277126 & -464976 &  337018 &  -88297 \\
   99213 & -464976 &  812586 & -611976 &  165153 \\
  -70237 &  337018 & -611976 &  485446 & -140251 \\
   18079 &  -88297 &  165153 & -140251 &   45316
\end{array}
\right],
\\
\mathbf{B}_2 &=& \scriptsize\tfrac{1}{10080}
\left[
\begin{array}{rrrrr}
   13816 &  -51001 &   67923 &  -38947 &   8209 \\
  -51001 &  209926 & -299076 &  179098 & -38947 \\
   67923 & -299076 &  462306 & -299076 &  67923 \\
  -38947 &  179098 & -299076 &  209926 & -51001 \\
    8209 &  -38947 &   67923 &  -51001 &  13816
\end{array}
\right],
\\
\mathbf{B}_d &=& \scriptsize \tfrac{1}{10080}
\left[
\begin{array}{rrrrr}
   942016 &  -2964751 &   3597813 & -2004757 &   429679 \\
 -2964751 &   9449926 & -11584896 &  6499258 & -1399537 \\
  3597813 & -11584896 &  14319786 & -8079576 &  1746873 \\
 -2004757 &   6499258 &  -8079576 &  4577926 &  -992851 \\
   429679 &  -1399537 &   1746873 &  -992851 &   215836
\end{array}
\right], 
\end{eqnarray*}
and $\mathbf{B}_{4-k} = \mathrm{rev}(\mathbf{B}_k), \; k = 0,1$.
The smoothness indicator for the downwind substencil $S_5$ is then defined by 
$$ \beta_5 = \frac{1}{6} \left( \sum_{k=0}^4 \beta_k + \beta_d \right). $$
Finally, the global smoothness indicator $\tau$ \eqref{eq:GSI} employs the anti-symmetric vector 
$$ \bfc_\tau = (-1,9,-36,84,-126,126,-84,36,-9,1). $$


\section*{Code and data availability}
The custom codes generated during the current study are available from the corresponding author upon reasonable request.
No datasets were generated or analyzed in the current study.

\section*{Conflict of interest}
We declare that all authors have no known competing financial interests or personal ties that could have an impact on the work reported in this paper.
\section*{Acknowledgments}
\label{sec:Ack}
Gu is supported by the POSTECH Basic Science Research Institute Fund, whose NRF grant number is RS-2021-NR060139.
Wang is supported in part by the National Natural Science Foundation of China (12301530), the Shandong Provincial Qingchuang Science and Technology Project, and the startup funding from Ocean University of China.
Jung is supported by the National Research Foundation of Korea (NRF) under grant number 2021R1A2C3009648, the POSTECH Basic Science Research Institute under the NRF grant number 2021R1A6A1A10042944, and partially by an NRF grant funded by the Korea government (MSIT) (RS-2023-00219980).
Don thanks the Hong Kong Research Grant Council (GRF/12301824, GRF/12300922).
Additionally, the authors, Wang and Don, extend their gratitude to POSTECH for providing a conducive environment for conducting this research during their visits.

\section*{Authorship contribution statement}
All the authors have contributed equally to this scholarly research and article.

\bibliographystyle{plain}
\bibliography{reference}

\end{document}